\RequirePackage{fix-cm}
\documentclass{svjour3}                     %

\smartqed  %

\usepackage{graphicx}

\usepackage[utf8]{inputenc} %
\usepackage[T1]{fontenc} %
\usepackage[english]{babel} %
\usepackage{graphicx,setspace,parskip,bm,amsmath} %
\usepackage{tabularx,multicol,multirow,threeparttablex} %
\usepackage{amssymb} %
\usepackage{footnpag,marginnote} %
\usepackage{enumerate} %
\usepackage[inline]{enumitem}  %
\usepackage{booktabs} %
\usepackage[font=small]{caption} %
\usepackage{subcaption} %
\usepackage{environ,tikz} %
\usetikzlibrary{mindmap}
\usepackage{mathptmx} %
\usepackage{microtype} %
\usepackage{makecell} %
\usepackage{setspace} %
\usepackage{xcolor} %
\usepackage{colortbl} %
\usepackage[top=2.5cm,bottom=2.5cm,left=2.5cm,right=2.5cm]{geometry} %
\usepackage{authblk} %
\usepackage{pdflscape} %
\usepackage{float} %
\usepackage[hidelinks,colorlinks]{hyperref} %
\usepackage{url} %
\usepackage{doi} %
\usepackage{csvsimple} %
\usepackage{siunitx} %
\usepackage{comment} %
\usepackage{appendix} %
\usepackage[pagewise]{lineno} %
\usepackage[ruled,vlined]{algorithm2e} %
\usepackage{stmaryrd} %
\usepackage{wasysym} %
\usepackage{etoolbox} %
\usepackage[nonumberlist,acronym]{glossaries} %
\usepackage{natbib}

\DeclareMathAlphabet{\mathcal}{OMS}{cmsy}{m}{n}

\graphicspath{{figs/}{figs/histograms/}{figs/regression/}{figs/violins/}{figs/results/}{figs/results/cartesian_mapping/}{figs/results/mapping/}{figs/results/mapping/hi_res/}{figs/results/mapping/lo_res/}{figs/results/mapping/GRATIS/}}

\hypersetup{
	linkcolor={cyan},
	citecolor={cyan},
	urlcolor={cyan}
}

\makeatletter
\newsavebox{\measure@tikzpicture}
\NewEnviron{scaletikzpicturetowidth}[1]{%
	\def\tikz@width{#1}%
	\def\tikzscale{1}\begin{lrbox}{\measure@tikzpicture}%
		\BODY
	\end{lrbox}%
	\pgfmathparse{#1/\wd\measure@tikzpicture}%
	\edef\tikzscale{\pgfmathresult}%
	\BODY
}
\makeatother

\tikzset{grow cyclic absolute/.style={%
		growth function=\tikzgrowcyclicabsolute}}
\def\tikzgrowcyclicabsolute{%
	\pgftransformshift{%
		\pgfpointpolar{\pgfkeysvalueof{/tikz/sibling angle}}%
		{\the\tikzleveldistance}}}

\def\anglemap{50}
\def\lengthmap{20}
\def\anglemapd{25}
\def\lengthmapd{7}

\DeclareSIUnit{\AU}{AU}
\DeclareSIUnit{\day}{d}
\DeclareSIUnit{\deg}{deg}
\DeclareSIUnit{\year}{yr}

\newcommand{\sgn}[1]{\mathrm{sgn}({#1})}

\newcommand{\Eq}[1]{Eq.~\eqref{#1}}
\newcommand{\Eqs}[1]{Eqs.~\eqref{#1}}
\newcommand{\Fig}[1]{Fig.~\ref{#1}}

\newcommand{\Tab}[1]{Table~\ref{#1}}
\newcommand{\Alg}[1]{Algorithm~\ref{#1}}
\newcommand{\Sec}[1]{Section~\ref{#1}}
\newcommand{\Subsec}[1]{Section~\ref{#1}}

\newcommand{\lastdate}{Feb 1, 2022}

\newcommand{\WSBset}[3]{$\mathcal{#1}_{#2}^{#3}$}
\newcommand{\Wset}[1]{\WSBset{W}{#1}{}}
\newcommand{\Xset}[1]{\WSBset{X}{#1}{}}
\newcommand{\Kset}[1]{\WSBset{K}{#1}{}}
\newcommand{\Mset}[1]{\WSBset{M}{#1}{}}
\newcommand{\Dset}[1]{\WSBset{D}{#1}{}}
\newcommand{\Cset}[1]{\WSBset{C}{-1}{#1}}
\newcommand{\DAIset}[1]{\WSBset{\hat{I}}{#1}{}}
\newcommand{\DAKset}[1]{\WSBset{\hat{K}}{#1}{}}
\newcommand{\DAMset}[1]{\WSBset{\hat{M}}{#1}{}}
\newcommand{\DAXset}[1]{\WSBset{\hat{X}}{#1}{}}
\newcommand{\DAWset}[1]{\WSBset{\hat{W}}{#1}{}}

\newcommand{\Sun}{\odot}

\newcommand{\textfnc}[1]{{\fontfamily{pcr}\selectfont #1}}

\newacronym{ABM}{ABM}{Adams-Bashforth-Moulton}
\newacronym{ADS}{ADS}{automatic domain splitting}
\newacronym{BC}{BC}{ballistic capture}
\newacronym{BCC}{BCC}{ballistic capture corridor}
\newacronym{BME}{BME}{body mean equator of date frame}
\newacronym{CK}{CK}{orientation or C-matrix kernel}
\newacronym{CR3BP}{CR3BP}{circular restricted three body problem}
\newacronym{DA}{DA}{differential algebra}
\newacronym{DACE}{DACE}{differential algebra core engine}
\newacronym{DAER}{DAER}{Department of Aerospace Science and Technology}
\newacronym{DOPRI8}{DOPRI8}{Dormand-Prince 8th-order embedded Runge-Kutta method}
\newacronym{DOP853}{DOP853}{Dormand-Prince of order 8(5,3) embedded Runge-Kutta method}
\newacronym{ECSS}{ECSS}{European Cooperation for Space Standardization}
\newacronym{EOM}{EoM}{equations of motion}
\newacronym{ERC}{ERC}{European Research Council}
\newacronym{FK}{FK}{frame kernel}
\newacronym{GMAT}{GMAT}{General Mission Analysis Tool}
\newacronym{GNC}{GNC}{guidance, navigation, and control}
\newacronym{GRATIS}{GRATIS}{GRAvity TIdal Slide}
\newacronym{IC}{IC}{initial condition}
\newacronym{ICRF}{ICRF}{International Celestial Reference Frame}
\newacronym{IK}{IK}{instrument kernel}
\newacronym{JPL}{JPL}{Jet Propulsion Laboratory}
\newacronym{LSK}{LSK}{leap seconds kernel}
\newacronym{NAIF}{NAIF}{Navigation and Ancillary Information Facility}
\newacronym{NASA}{NASA}{National Aeronautics and Space Administration}
\newacronym{NEA}{NEA}{near-Earth asteroid}
\newacronym{NEO}{NEO}{ear-Earth object}
\newacronym{NSG}{NSG}{non-spherical gravity}
\newacronym{PECE}{PECE}{predictor-corrector}
\newacronym{PCK}{PCK}{planetary constants kernel}
\newacronym{PDS}{PDS}{Planetary Data System}
\newacronym{RAAN}{RAAN}{right ascension of the ascending node}
\newacronym{RHS}{RHS}{right-hand side}
\newacronym{RK}{RK}{Runge-Kutta}
\newacronym{RSS}{RSS}{root sum square}
\newacronym{RTN}{RTN@$t_{i}$}{radial-tangential-normal of date frame}
\newacronym{RPF}{RPF}{roto-pulsating frame}
\newacronym{SCLK}{SCLK}{spacecraft clock kernel}
\newacronym{SPK}{SPK}{spacecraft ephemeris kernel}
\newacronym{SRP}{SRP}{solar radiation pressure}
\newacronym{TSAE}{TSAE}{Toulouse graduate School of Aerospace Engineering}
\newacronym{UOM}{uom}{unit of measurement}
\newacronym{VSVO}{VSVO}{variable-step variable-order}
\newacronym{VV}{V\&V}{verification and validation}
\newacronym{WSB}{WSB}{weak stability boundary}

 \journalname{Celestial Mechanics and Dynamical Astronomy}
\begin{document}
\title{Stable sets mapping with Taylor differential algebra with application to ballistic capture orbits around Mars}

\author{T. Caleb \and
G. Merisio \and
P. di Lizia \and
F. Topputo
}

\institute{Thomas Caleb (corresponding author) \at
              ISAE-SUPAERO, 10 Avenue Édouard Belin, 54032 - 31055, Toulouse, France\\
              OrcID: 0000-0002-4027-3340\\
              \email{thomas.caleb@student.isae-supaero.fr}
           \and
           Gianmario Merisio \at
              Department of Aerospace Science and Technology, Politecnico di Milano, Via La Masa 34, 20156, Milano, Italy\\
              OrcID: 0000-0001-8806-7952\\
              \email{gianmario.merisio@polimi.it}
           \and
           Pierluigi di Lizia  \at
           	  Department of Aerospace Science and Technology, Politecnico di Milano, Via La Masa 34, 20156, Milano, Italy\\
           	    OrcID: 0000-0003-1692-3929\\
              \email{pierluigi.dilizia@polimi.it}
           \and
           Francesco Topputo  \at
              Department of Aerospace Science and Technology, Politecnico di Milano, Via La Masa 34, 20156, Milano, Italy\\
              OrcID: 0000-0002-5369-6887\\
              \email{francesco.topputo@polimi.it}
}

\date{Received: date / Accepted: date}

\maketitle              %
\begin{abstract}
Ballistic capture orbits offer safer Mars injection at longer transfer time. However, the search for such an extremely rare event is a computationally-intensive process. Indeed, it requires the propagation of a grid sampling the whole search space. This work proposes a novel ballistic capture search algorithm based on Taylor differential algebra propagation. This algorithm provides a continuous description of the search space compared to classical grid sampling research and focuses on areas where the nonlinearities are the largest.Macroscopic analyses have been carried out to obtain cartography of large sets of solutions. Two criteria, named consistency and quality, are defined to assess this new algorithm and to compare its performances with classical grid sampling of the search space around Mars. Results show that differential algebra mapping works on large search spaces, and automatic domain splitting captures the dynamical variations on the whole domain successfully. The consistency criterion shows that more than 87\% of the search space is guaranteed as accurate, with the quality criterion kept over 80\%.

\end{abstract}

\section{Introduction} \label{sec:intro}

\Gls*{BC} allows a spacecraft to approach a planet and enter a temporary orbit about it without requiring maneuvers in between. As part of the low-energy transfers, it is a valuable alternative to Keplerian approaches. Exploiting \gls*{BC} grants several benefits in terms of both cost reduction \citep{belbruno1993sun} and mission versatility \citep{belbruno2000calculation,topputo2015earth}, in general at the cost of longer transfer times \citep{circi2001dynamics,ivashkin2002trajectories}. In the past, the \gls*{BC} mechanism was used to rescue Hiten \citep{belbruno1990ballistic}, and to design insertion trajectories in lunar missions like SMART-1 \citep{racca2002smart} and GRAIL \citep{chung2010trans}. In the near future, BepiColombo will exploit \gls*{BC} orbits to be weakly captured by Mercury \citep{benkhoff2021bepicolombo,schuster2014influence}. \gls*{BC} is an event occurring in extremely rare occasions and requires acquiring a proper state (position and velocity) far away from the target planet \citep{topputo2015earth}. In fact, massive numerical simulations are required to find the specific conditions that support capture \citep{topputo2009computation} and only approximately 1 out of \num{10000} states lead to capture \citep{luo2015analysis}. In a first effort to reduce the computational burden, the variational theory for Lagrangian coherent structures \citep{haller2011variational} was recently applied to find \gls*{BC} opportunities more efficiently \citep{manzi2021flow}.

\Gls*{DA} propagation is a worthy candidate to reduce the computational burden for the search of \gls*{BC} trajectories. It consists of propagating \glspl*{IC}, not as a single point but as an interval around an \gls*{IC}. Thanks to the Taylor expansions of the flow, the state of any point in the represented interval can be determined through convenient polynomial evaluations \citep{berz1992high,berz1999modern}. This gain in efficiency comes at the cost of a loss of accuracy due to the finite Taylor expansion. The \gls*{ADS} algorithm allows to represent large domains accurately when increasing the order of the polynomial expansions fails to do so \citep{wittig2015propagation}.  Indeed, on poorly-defined domains, raising the order will increase the approximation error on the edges while the error will decrease on already well-defined areas. \Gls*{ADS} splits the initial domain into smaller sub-domains to reduce the approximation error when it grows above a given tolerance.

\gls*{DA} propagation is increasingly used in astrodynamics. Indeed, it provides high performances in uncertainty propagation in the two-body dynamics \citep{valli2013nonlinear}, even for highly nonlinear dynamics with large uncertainties when exploiting \gls*{ADS}, as in the case of Apophis \citep{wittig2015propagation} or Apollo LM-10 also known as Snoopy \citep{caleb2021can}. Other applications occur in orbital mechanics, such as propagation of probability density functions \citep{wittig2017longterm}, maximum a posteriori estimation \citep{servadio2022maximum}, orbit determination \citep{pirovano2021differential,servadio2021differential}, and generation and study of orbit families in the \gls*{CR3BP} \citep{dilizia2008application,baresi2021highorder}.

\gls*{DA} propagation allows avoiding intensive grid sampling implied by point-wise research of conventional algorithms for designing \gls*{BC} trajectories \citep{hyeraci2010method,luo2014constructing}, as it offers a continuous description of the whole search-space. In addition, the computation of a Taylor expansion provides information such as the partial derivatives of the flow up to an arbitrary order \citep{wittig2017longterm}. Furthermore, the polynomial maps can be manipulated to impose constraints on the flown trajectories \citep{berz1999modern, dilizia2008application}.

The goal of this work is to use \gls*{DA} mapping to carry out macroscopic analyses of the phase space about Mars to find \gls*{BC} trajectories. Hence, the adaptation of the definition of \gls*{BC} for \gls*{DA}, followed by the definition of two criteria, named consistency and quality, to assess the performances of \gls*{DA}-based mapping of \gls*{BC} compared to point-wise mapping. Cartography of large \gls*{BC} sets about Mars are computed using \gls*{DA} mapping. Mars is chosen without loss of generality due to its relevance in the long-term exploration. The work proposes an alternative classification algorithm to sort the sub-domains produced by the \gls*{ADS} algorithm in newly defined capture sets. While \gls*{DA} mapping allows performing good macroscopic cartography of the search space, it does not allow to fully replace point-wise mapping. However, \gls*{DA} mapping has the advantage of being continuous, so meaning that the behavior of any point in the search space is defined. On the contrary, a point-wise mapping delivers precise knowledge on the discrete set of propagated \glspl*{IC} only. As a consequence, the behavior of the continuous space between two points is unknown and can be challenging to interpolate due to the nonlinear dynamics. 

The remainder of the paper is organized as follows. In \Sec{sec:background}, the dynamical model employed is introduced, as well as the \gls*{WSB} concept, the \gls*{BC} mechanism, and the \gls*{DA} propagation. Then, the description of the characterization process, and the mapping-assessment methodology follow in \Sec{sec:method}. Results are presented and discussed in \Sec{sec:results}. Eventually, conclusions are drawn in \Sec{sec:conclusion} together with the presentation of future work.

\section{Background} \label{sec:background}

Details about dynamical model, \gls*{WSB}, \gls*{BC} phenomenon, \gls*{DA} propagation, and \gls*{ADS} algorithms are herewith presented.

\subsection{Dynamical model}

According to the nomenclature introduced in \citet{luo2014constructing}, a \textit{target} (also referred to as \textit{central body}) and a \textit{primary} are defined. The target being the body around which the motion of the spacecraft is studied (Mars in this work), and the primary being the body around which the target revolves (the Sun). Target and primary masses are $m_{t}$ and $m_{p}$, respectively.

\subsubsection{Reference frames} \label{sec:reference_frame}
In this work, the following reference frames are used: J2000, and \acrshort*{RTN}.

\paragraph{J2000.}
Defined on the Earth's mean equator and equinox, the J2000 is an inertial frame determined from observations of planetary motions which was realized to coincide almost exactly with the \gls*{ICRF} \citep{archinal2011report}. \Gls*{EOM} are integrated in this reference frame.

\paragraph{RTN@$t_{i}$.}
The \gls*{RTN} is an inertial frame frozen at a prescribed epoch $t_{i}$. The frame is centered at the target. The $x$-axis is aligned with the primary--secondary direction, the $z$-axis is normal to the primary--secondary plane in the direction of their angular momentum, and the $y$-axis completes the dextral orthonormal triad. \glspl*{IC} are defined in this frame \citep{luo2015analysis}.

\subsubsection{Ephemerides}
Precise states of the Sun and the major planets are retrieved from the \gls*{JPL}'s planetary ephemerides \textfnc{de440s.bsp}\footnote{Data publicly available at: \url{https://naif.jpl.nasa.gov/pub/naif/generic_kernels/spk/planets/de440s.bsp} [retrieved \lastdate].} (or DE440s) \citep{park2021jpl}. Additionally, the ephemerides \textfnc{mars097.bsp} of Mars (the target) and its moons are employed\footnote{\url{~/spk/satellites/mars097.bsp} [retrieved \lastdate].}. The following generic \gls*{LSK} and \gls*{PCK} are used: \textfnc{naif0012.tls}, \textfnc{pck00010.tpc}, and \textfnc{gm\_de431.tpc}\footnote{Data publicly available at: \url{https://naif.jpl.nasa.gov/pub/naif/generic_kernels/lsk/naif0012.tls}, \url{~/generic_kernels/pck/pck00010.tpc}, and \url{~/generic_kernels/pck/gm_de431.tpc} [retrieved \lastdate]. The \textfnc{gm\_de431.tpc} \gls*{PCK} kernel is used because the new version consistent with the ephemerides DE440s has not been released yet.}.

\subsubsection{Equations of motion}
The \gls*{EOM} used are those of the restricted $N$-body problem. The gravitational attractions of the Sun, Mercury, Venus, Earth (B\footnote{Here B stands for barycenter.}), Mars (central body), Jupiter (B), Saturn (B), Phobos, and Deimos are considered. Additionally, \gls*{SRP} is also included and implemented as a \emph{cannonball} or \emph{spherical} model \citep{scheeres2011dynamics}. The assumed spacecraft specifications needed to evaluate the \gls*{SRP} perturbation are collected in \Tab{tab:sc-specs}. They are compatible with the specifications of a 12U deep-space CubeSat \citep{topputo2021envelop}.

\begin{table}[tbp]
	\centering
	\caption{Assumed spacecraft specifications.}
	\begin{tabular}{lll}
		\toprule
		\textbf{Specification} & \textbf{Symbol} & \textbf{Value} \\
		\midrule
		Mass & $ m $ & \SI{24}{\kilo\gram} \\
		\gls*{SRP} area & $ A $ & \SI{0.32}{\meter\squared} \\
		Coefficient of reflectivity & $ C_{r} $ & \SI{1.3}{} \\
		\bottomrule
	\end{tabular}
	\label{tab:sc-specs}
\end{table}  

The \gls*{EOM}, written in a non-rotating Mars-centered reference frame are \citep{luo2014constructing,merisio2021characterization}
\begin{equation}
	\ddot{\mathbf{r}} = - \frac{\mu_{t}}{r^{3}} \mathbf{r} 
	- \sum\limits_{i \in \mathbb{P}} \mu_{i} \left( \frac{\mathbf{r}_{i}}{r_{i}^{3}} + \frac{\mathbf{r} - \mathbf{r}_{i}}{\left\lVert \mathbf{r} - \mathbf{r}_{i} \right\rVert^{3}} \right)
	+ \frac{Q A}{m} \frac{\mathbf{r} - \mathbf{r}_{\Sun}}{\left\lVert \mathbf{r} - \mathbf{r}_{\Sun} \right\rVert^{3}}
	\label{eq:eom}
\end{equation}
where $ \mu_{t} $ is the gravitational parameter of the target body; $ \mathbf{r} $ is the position vector of the spacecraft with respect to the target and $ r $ is its magnitude; $ \mathbb{P} $ is a set of $ N-2 $ indexes each referring to a perturbing body; $ \mu_{i} $ and $ \mathbf{r}_{i} $ are the gravitational parameter and position vector with respect to the target of the $i$-th body, respectively; $ A $ is the Sun-projected area on the spacecraft for \gls*{SRP} evaluation; $ m $ is the spacecraft mass; $ \mathbf{r}_{\Sun} $ is the position vector of the Sun with respect to the target. Lastly, $ Q $ is equal to
\begin{equation}
	Q = \frac{L C_{r}}{4 \pi c} 
	\label{eq:solarq}
\end{equation}
where $ C_{r} $ is the spacecraft coefficient of reflectivity, $ c = \SI{299792458}{\meter\per\second} $ taken from SPICE \citep{acton1996ancillary,acton2018look} is the speed of light in vacuum, and $ L = S_{\Sun} 4 \pi d_{\mathrm{AU}}^{2} $ is the luminosity of the Sun. The latter is computed from the solar constant\footnote{\url{https://extapps.ksc.nasa.gov/Reliability/Documents/Preferred_Practices/2301.pdf} [last accessed \lastdate].} $ S_{\Sun} = \SI{1367.5}{\watt\per\meter\squared} $  evaluated at $ d_{\mathrm{AU}} = \SI{149597870613.6889}{\meter} $ corresponding to \SI{1}{\AU} according to SPICE \citep{acton1996ancillary,acton2018look}.

\subsubsection{Numerical integration of \gls*{EOM}}
The \gls*{EOM} in \Eq{eq:eom} are integrated in their nondimensional form to avoid ill-conditioning \citep{luo2014constructing}. Nondimensionalization units are reported in \Tab{tab:units}. For point-wise simulations, the numerical integration is carried with the \gls*{DOPRI8} propagation scheme \citep{montenbruck2000satellite}. It is an adaptive step, 8th-order \gls*{RK} integrator with 7th-order error control, the coefficients were derived by Prince and Dormand \citep{prince1981high}. As for the \gls*{DA} simulations, the propagation scheme is \gls*{DOP853} \citep{hairer1993solving}, which is of the same 8th-order \gls*{RK} integrator family. However, the error control is performed with a 3rd-order and a 5th-order estimation. The dynamics are propagated setting the relative tolerance to $ 10^{-12} $ \citep{luo2014constructing}. 

\begin{table}[tbp]
	\centering
	\caption{Nondimensionalization units.}
	\begin{tabular}{llll}
		\toprule
		\textbf{Unit} & \textbf{Symbol} & \textbf{Value} & \textbf{Comment} \\
		\midrule
		Gravity parameter & $ \mathrm{MU} $ & \SI{42828.376}{\kilo\meter\cubed\per\second\squared} & Mars' gravity parameter $ \mu_{t} $ \\
		Length & $ \mathrm{LU} $ & \SI{3396.0000}{\kilo\meter} & Mars' radius $ R_{\mars} $ \\ 
		Time & $ \mathrm{TU} $ & \SI{956.28142}{\second} & $ (\mathrm{LU}^3 / \mathrm{MU})^{0.5} $ \\
		Velocity & $ \mathrm{VU} $ & \SI{3.5512558}{\kilo\meter\per\second} & $ \mathrm{LU} / \mathrm{TU} $ \\
		\bottomrule
	\end{tabular}
	\label{tab:units}
\end{table}

\subsection{Weak stability boundary and ballistic capture mechanism} \label{sec:wsb}

Over the years, the \gls*{WSB} was defined in many different ways. It was initially identified as a fuzzy boundary region placed at approximately \SI{1.5e6}{\kilo\meter} from the Earth in the Sun--Earth direction \citep{belbruno1987lunar,belbruno1990ballistic}. An algorithmic definition followed in \citet{belbruno2004capture}, later extended in \citet{garcia2007note}, \citet{topputo2009computation}, and \citet{silva2012applicability}. Then, the \gls*{WSB} was interpreted as the intersection of three sub-sets of the phase space \citep{topputo2008resonant,belbruno2008resonance}. The \gls*{WSB} concept being closely connected to \gls*{BC} \citep{belbruno2004capture}, a formal definition and a methodology for its derivation from weakly stable and unstable sets were finally proposed in \citet{hyeraci2010method}. To date, despite the effort put in numerous works \citep{garcia2007note,topputo2008resonant,belbruno2008resonance,belbruno2010weak,belbruno2013geometry}, both \gls*{WSB} and \gls*{BC} are still not completely understood. Nonetheless, a connection between celestial and quantum mechanics was recently found exploiting the \gls*{WSB} \citep{belbruno2020relation}, providing a fresh perspective to tackle the problem.

\Gls*{BC} orbits are characterized by \glspl*{IC} escaping the target when integrated backward and performing $n$ revolutions about it when propagated forward, neither impacting or escaping the target. In forward time, particles flying on \gls*{BC} orbits approach the target coming from outside its sphere of influence and remain temporarily captured about it. After a certain time, the particle escapes if an energy dissipation mechanism does not take place to make the capture permanent. To dissipate energy either a breaking maneuver or the target's atmosphere (if available) could be used \citep{luo2021mars}. In this work, \gls*{BC} sets for comparison purposes are derived propagating the \gls*{EOM} in \Eq{eq:eom} and following the procedure in \citet{luo2014constructing}.

When searching for \gls*{BC} opportunities, most of the trajectories found are spurious solutions which are typically not useful for mission design purposes \citep{luo2014constructing}. Useful solutions are detected exploiting the regularity index\footnote{In previous works this was referred to as stability index \citep{luo2014constructing,luo2015analysis,luo2017capability}. However, in \citet{deitos2018survey}, the adjustment from \emph{stability} to \emph{regularity} index was proposed to avoid misunderstandings with the periodic orbit stability index. The same nomenclature introduced in \citet{deitos2018survey} is used in this work.} $ S $ and regularity coefficient $ \Delta S_{\%} $ \citep{deitos2018survey}. The aim is seeking for ideal orbits that presents regular post-capture legs resulting in $n$ revolutions about the target which are similar in orientation and shape. Numerical experiments showed that high-quality post-capture orbits are identified by small regularity index and coefficient \citep{deitos2018survey,luo2017capability,luo2015analysis,luo2014constructing}. If the regularity index and coefficient are indicators used to qualitatively judge post-capture legs, capture occurrence is quantitatively measured through the capture ratio $ \mathcal{R}_{\mathcal{C}} $ \citep{luo2015analysis}. Typically, search spaces characterized by larger capture ratio are desirable when looking for \gls*{BC} orbits.

Effects on \gls*{BC} by gravitational attractions of many bodies besides the primaries and \gls*{SRP} have been investigated in previous works \citep{merisio2021characterization,aguiar2018technique,luo2017capability,luo2015analysis}, with \citet{merisio2021characterization}, and \citet{aguiar2018technique} also considering the target non-spherical gravity perturbations. Compared to the restricted three-body problem, the $N$-body model is more adequate for constructing ballistic capture orbits as proved in \citet{luo2015analysis}, with a particle being more easily captured when considering additional gravitational attractions. The latter was later confirmed in \citet{luo2017capability}, where the presence of moons is demonstrated to increase the capture ratio $\mathcal{R}_{\mathcal{C}}$ and improves the regularity of post-capture orbits while accommodating larger pre-capture energies. Similarly, \gls*{SRP} increases the chances of being temporarily captured about the target and regularizes the post-capture legs \citep{merisio2021characterization,aguiar2018technique}. Even if not duly discussed, the aforementioned remarks are observed also in the results of this work. Overall, supplementary terms in the dynamics seem to favor the manifestation of the \gls*{BC} phenomenon.

\subsubsection{Definitions of particle stability and sub-sets} \label{sec:set_defs}
A particle stability is inferred using a plane in the three-dimensional physical space \citep{belbruno1993sun}, according to the spatial stability definition provided in \citet{luo2014constructing}. The following indications are used to classify stability, see \citet{luo2014constructing} for more details:
\begin{enumerate*}[label=\arabic*)]
	\item a particle completes a revolution around the target according to \textit{Remark 1} and Eq.~(5) in \citet{luo2014constructing};
	\item a particle escapes from the target according to \textit{Remark 2} and Eq.~(6) in \citet{luo2014constructing};
	\item a particle impacts with the target according to \textit{Remark 3} and Eq.~(7) in \citet{luo2014constructing}.
\end{enumerate*}
Consistent variants of Eq.~(7) in \citet{luo2014constructing} can be derived to locate impacts with target's moons, if present.

\begin{figure}[tbp]
	\begin{center}
		\begin{scaletikzpicturetowidth}{\textwidth}
			\begin{tikzpicture}[scale=\tikzscale]
			\node[rotate=90, align=left, text width=20, anchor=south west] at (-0.05,0) {Backward};
			\node[rotate=90, align=left, text width=20, anchor=north west] at (0.05,0) {Forward};
			\draw[line width=0.5pt] (0,0) -- (0,1);
			
			\node  at (0.1,0.9) {$ \Omega $}[grow cyclic absolute]
			child[sibling angle=0,level distance=\lengthmap]{node{\Wset{1}}
				child[sibling angle=0,level distance=\lengthmap]{node{\Wset{2}}
					child[dashed,sibling angle=0,level distance=\lengthmapd]{node{}}
					child[dashed,sibling angle=-atan(tan(\anglemapd)*1/4),level distance=\lengthmapd]{node{}}
					child[dashed,sibling angle=-atan(tan(\anglemapd)*2/4),level distance=\lengthmapd]{node{}}
					child[dashed,sibling angle=-atan(tan(\anglemapd)*3/4),level distance=\lengthmapd]{node{}}
					child[dashed,sibling angle=-\anglemapd,level distance=\lengthmapd]{node{}}
				}
				child[sibling angle=-atan(tan(\anglemap)*1/4),level distance=sqrt(\lengthmap^2+(\lengthmap*tan(\anglemap)*1/4)^2)]{node{\Xset{2}}}
				child[sibling angle=-atan(tan(\anglemap)*2/4),level distance=sqrt(\lengthmap^2+(\lengthmap*tan(\anglemap)*2/4)^2)]{node{\Kset{2}}}
				child[sibling angle=-atan(tan(\anglemap)*3/4),level distance=sqrt(\lengthmap^2+(\lengthmap*tan(\anglemap)*3/4)^2)]{node{\Mset{2}}}
				child[sibling angle=-\anglemap,level distance=\lengthmap/cos(\anglemap)]{node{\Dset{2}}}
			}
			child[sibling angle=-atan(tan(\anglemap)*1/4),level distance=sqrt(\lengthmap^2+(\lengthmap*tan(\anglemap)*1/4)^2)]{node{\Xset{1}}}
			child[sibling angle=-atan(tan(\anglemap)*2/4),level distance=sqrt(\lengthmap^2+(\lengthmap*tan(\anglemap)*2/4)^2)]{node{\Kset{1}}}
			child[sibling angle=-atan(tan(\anglemap)*3/4),level distance=sqrt(\lengthmap^2+(\lengthmap*tan(\anglemap)*3/4)^2)]{node{\Mset{1}}}
			child[sibling angle=-\anglemap,level distance=\lengthmap/cos(\anglemap)]{node{\Dset{1}}};
			
			\node at (-0.1,0.9) {$ \Omega $}[grow cyclic absolute] 
			child[sibling angle=180,level distance=\lengthmap]{node{\Wset{-1}}
				child[sibling angle=180,level distance=\lengthmap]{node{\Wset{-2}}
					child[dashed,sibling angle=180,level distance=1.09*\lengthmapd]{node{}}
					child[dashed,sibling angle=180+atan(tan(\anglemapd)*1/4),level distance=1.09*\lengthmapd]{node{}}
					child[dashed,sibling angle=180+atan(tan(\anglemapd)*2/4),level distance=1.09*\lengthmapd]{node{}}
					child[dashed,sibling angle=180+atan(tan(\anglemapd)*3/4),level distance=1.09*\lengthmapd]{node{}}
					child[dashed,sibling angle=180+\anglemapd,level distance=1.1*\lengthmapd]{node{}}
				}
				child[sibling angle=180+atan(tan(\anglemap)*1/4),level distance=sqrt(\lengthmap^2+(\lengthmap*tan(\anglemap)*1/4)^2)]{node{\Xset{-2}}}
				child[sibling angle=180+atan(tan(\anglemap)*2/4),level distance=sqrt(\lengthmap^2+(\lengthmap*tan(\anglemap)*2/4)^2)]{node{\Kset{-2}}}
				child[sibling angle=180+atan(tan(\anglemap)*3/4),level distance=sqrt(\lengthmap^2+(\lengthmap*tan(\anglemap)*3/4)^2)]{node{\Mset{-2}}}
				child[sibling angle=180+\anglemap,level distance=\lengthmap/cos(\anglemap)]{node{\Dset{-2}}}
			}
			child[sibling angle=180+atan(tan(\anglemap)*1/4),level distance=sqrt(\lengthmap^2+(\lengthmap*tan(\anglemap)*1/4)^2)]{node{\Xset{-1}}}
			child[sibling angle=180+atan(tan(\anglemap)*2/4),level distance=sqrt(\lengthmap^2+(\lengthmap*tan(\anglemap)*2/4)^2)]{node{\Kset{-1}}}
			child[sibling angle=180+atan(tan(\anglemap)*3/4),level distance=sqrt(\lengthmap^2+(\lengthmap*tan(\anglemap)*3/4)^2)]{node{\Mset{-1}}}
			child[sibling angle=180+\anglemap,level distance=\lengthmap/cos(\anglemap)]{node{\Dset{-1}}};
			\end{tikzpicture}
		\end{scaletikzpicturetowidth}
	\end{center}
	\caption{Sets relations according to the algorithmic definition of WSB.}
	\label{fig:graph}
\end{figure}
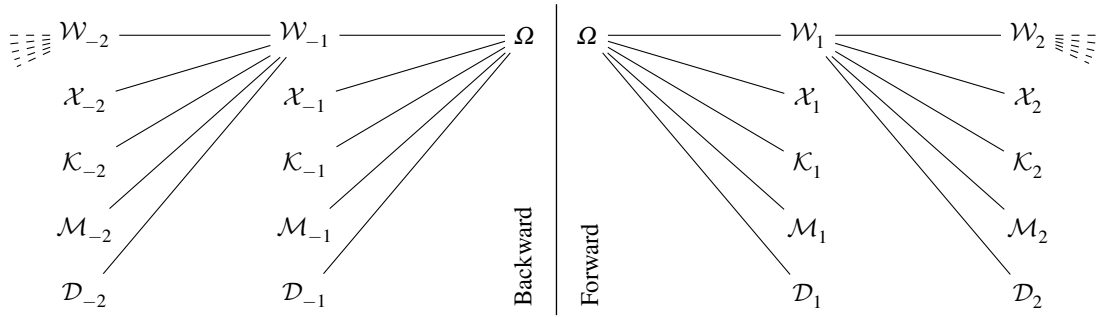

Based on its dynamical behavior, a propagated trajectory is said to be:
\begin{enumerate*}[label=\roman*)]
	\item \textit{weakly stable} (sub-set \Wset{i}) if the particle performs $i$ complete revolutions around the target without escaping or impacting with it or its moons;
	\item \textit{unstable} (sub-set \Xset{i}) if the particle escapes from the target before completing the $i$-th revolution;
	\item \textit{target--crash} (sub-set \Kset{i}) if the particle impacts with the target before completing the $i$-th revolution;
	\item \textit{moon--crash} (sub-set \Mset{i}) if the particle impacts with one of the target's moons before completing the $i$-th revolution;
	\item \textit{acrobatic} (sub-set \Dset{i}) if none of the previous conditions occurs within the integration time span.
\end{enumerate*}
Conditions ii)-v) apply after the particle performs $(i-1)$ revolutions around the target. The sub-sets are defined for $ i \in \mathbb{Z} \textbackslash \{0\} $, where the sign of $i$ informs on the propagation direction. When $ i > 0 $ ($ i < 0 $) the \gls*{IC} is propagated forward (backward) in time. The overall domain, union of all sub-sets, is defined as $\Omega$. A graph clarifying the relations between sub-sets is shown in \Fig{fig:graph}. A capture set is defined as $ \mathcal{C}_{-1}^{n} := \mathcal{W}_{n} \cap \mathcal{X}_{-1} $. Therefore, it is the intersection between the stable set in forward time \Wset{n} and the unstable set in backward time \Xset{-1} \citep{luo2014constructing}.

\subsection{Differential algebra} \label{sec:da}

\gls*{DA} propagation consists of assimilating a function $f$ of $v$ variables, contained in $C^{k+1}$, with $T_f^{(k)}$ the Taylor expansion of $f$ at order $k$ \citep{berz1999modern}. The computation of such polynomials can be performed efficiently, and provides a representation of the function $f$ on all of its domain. Moreover, the \gls*{DA} structure ensures that algebraic and functional operations are well-defined, particularly for the numerical solving of ordinary differential equations \citep{berz1992high}.

The main advantage of this method is that the polynomial map only needs to be computed once, and then it is evaluated in an arbitrarily large number of points. In other words, to perform a Monte-Carlo estimation with a sample of size $S$, only one computation of the map is needed, followed by $S$ polynomial evaluations, while classic Monte-Carlo requires $S$ propagations, see \citet{armellin2010asteroid}. The \gls*{DA} engine used in this work is the \gls*{DACE}\footnote{Library available at: \url{https://github.com/dacelib/dace} [last accessed \lastdate].}, implemented by Politecnico di Milano \citep{rasotto2016differential,massari2018differential}.

\subsubsection{Automatic domain splitting} \label{sec:ads}

DA propagation allows reducing the approximation error by increasing the order $k$ of the polynomial mapping. Nevertheless, increasing the order leads to an important growth of the computational time. Thus, the introduction of \gls*{ADS}, see \citet{wittig2015propagation}, and \citet{pirovano2021differential} for more details.

The farther a point is from the constant part of the Taylor expansion, the higher the loss of accuracy. The \gls*{ADS} algorithm controls this error by dividing the domain into halves until the accuracy is below a predetermined tolerance. Therefore, each new sub-domain is represented by its polynomial map causing a much smaller approximation error, at a controlled computational cost \citep{wittig2015propagation}. The \gls*{ADS} creates a division of the initial domain separating areas with different behaviors from one another. The \gls*{ADS} algorithm requires the passing of two parameters, in addition to the order of the \gls*{DA}: the tolerance, and the maximum number of splits allowed. The former being the approximation error threshold before the splitting of a domain occurs, while the latter being the number of times the \gls*{ADS} routine can be applied to a sub-domain.

\section{Methodology} \label{sec:method}

To map \gls*{BC} sets on the search space using \gls*{DA} propagation, a novel classification algorithm is devised and herewith presented. To be successful, the algorithm requires an accurate revolution period estimation. Consistency and quality criteria are defined to evaluate performances of the resulting mapping. Finally, fresh representation methods to visualize properly the large resulting amount of data are exposed.

The search space is chosen to maximize the capture ratio $ \mathcal{R}_{c} $ based on the analysis reported in \citet{luo2015analysis}. It is defined in the Mars-centered RTN reference frame at capture epoch $ t_{0} $ on December 9, 2023, at 00:45:18.363 (UTC). At that epoch, Mars's true anomaly with respect to the Sun is equal to \SI{270}{\deg}, maximizing $ \mathcal{R}_{c} $ \citep{luo2015analysis}. The selected plane is defined by inclination $ i = \SI{0.6283}{\radian} $, and right ascension of the ascending node $ \mathrm{RAAN} = \SI{0.6283}{\radian}$. That because, according to Fig. 10 in \citet{luo2015analysis}, such values maximize the capture ratio for Mars. Sought trajectories have osculating eccentricity $ e = 0.99 $ \citep{topputo2015earth}, and mean anomaly $ M = \SI{0}{\radian} $ at the initial epoch $ t_{0} $. If $ R_{\mars} $ is the radius of Mars in \si{\kilo\meter}, then the search space on the plane defined above is a circular crown centered at Mars, from radius $ R_{\mars} + \SI{100}{\kilo\meter} $ up to radius $ 5 R_{\mars} $. Hence, 
\begin{equation}
    \label{eq:search_space}
    \left(r_p, \omega\right) \in \left[R_{\mars} + \text{ \SI{100}{\kilo\meter}}, 5 \cdot R_{\mars}\right] \times \left(-\pi,\pi\right]
\end{equation}
with $r_p$ the radius of the periapsis, and $\omega$ the argument of the periapsis.

The main difference between the point-wise mapping of sets performed by \gls*{GRATIS} \citep{topputo2018trajectory} and the \gls*{DA} mapping presented in this work is that the \gls*{DA} propagator does not allow to count revolutions. It means the sub-domains cannot be classified in the sub-sets defined in \Subsec{sec:set_defs}. New sets are defined to solve this problem. Furthermore, instead of tracking revolutions geometrically, as in \citet{luo2014constructing}, the proposed novel classification algorithm counts revolution periods. Finally, a bridge between the two mapping methods is established.

\subsection{Redefinition of sub-sets and classification algorithm} \label{sec:hatted_sets}
There is no direct method to use the definitions of \Subsec{sec:set_defs} with \gls*{DA}. This is because accurately counting the revolutions of a continuous set of particles around the target cannot be performed as in \textit{Remark 1} and Eq.~(5) in \citet{luo2014constructing}. Instead of using this condition to count completed revolutions, the number of period elapsed since epoch is used. Consequently, \Dset{i} cannot be defined when using \gls*{DA} propagation. Moreover, if a subdomain reached the minimum size allowed by the \gls*{ADS} algorithm, and it tries to split again, the accuracy of the mapping of this sub-domain cannot be guaranteed. Thus, these inconsistent sub-domains need to be ruled out from the sub-sets, so that only the consistent ones are retained.

\begin{figure}[tbp]
	\begin{center}
		\begin{scaletikzpicturetowidth}{\textwidth}
			\begin{tikzpicture}[scale=\tikzscale]
			\node[rotate=90, align=left, text width=20, anchor=south west] at (-0.05,0) {Backward};
			\node[rotate=90, align=left, text width=20, anchor=north west] at (0.05,0) {Forward};
			\draw[line width=0.5pt] (0,0) -- (0,1);
			
			\node  at (0.1,0.9) {$ \Omega $}[grow cyclic absolute]
			child[sibling angle=0,level distance=\lengthmap]{node{\DAWset{1}}
				child[sibling angle=0,level distance=\lengthmap]{node{\DAWset{2}}
					child[dashed,sibling angle=0,level distance=1.1*\lengthmapd]{node{}}
					child[dashed,sibling angle=-atan(tan(\anglemapd)*1/4),level distance=1.1*\lengthmapd]{node{}}
					child[dashed,sibling angle=-atan(tan(\anglemapd)*2/4),level distance=1.1*\lengthmapd]{node{}}
					child[dashed,sibling angle=-atan(tan(\anglemapd)*3/4),level distance=1.1*\lengthmapd]{node{}}
					child[dashed,sibling angle=-\anglemapd,level distance=1.1*\lengthmapd]{node{}}
				}
				child[sibling angle=-atan(tan(\anglemap)*1/4),level distance=sqrt(\lengthmap^2+(\lengthmap*tan(\anglemap)*1/4)^2)]{node{\DAXset{2}}}
				child[sibling angle=-atan(tan(\anglemap)*2/4),level distance=sqrt(\lengthmap^2+(\lengthmap*tan(\anglemap)*2/4)^2)]{node{\DAKset{2}}}
				child[sibling angle=-atan(tan(\anglemap)*3/4),level distance=sqrt(\lengthmap^2+(\lengthmap*tan(\anglemap)*3/4)^2)]{node{\DAMset{2}}}
				child[sibling angle=-\anglemap,level distance=\lengthmap/cos(\anglemap)]{node{\DAIset{2}}}
			}
			child[sibling angle=-atan(tan(\anglemap)*1/4),level distance=sqrt(\lengthmap^2+(\lengthmap*tan(\anglemap)*1/4)^2)]{node{\DAXset{1}}}
			child[sibling angle=-atan(tan(\anglemap)*2/4),level distance=sqrt(\lengthmap^2+(\lengthmap*tan(\anglemap)*2/4)^2)]{node{\DAKset{1}}}
			child[sibling angle=-atan(tan(\anglemap)*3/4),level distance=sqrt(\lengthmap^2+(\lengthmap*tan(\anglemap)*3/4)^2)]{node{\DAMset{1}}}
			child[sibling angle=-\anglemap,level distance=\lengthmap/cos(\anglemap)]{node{\DAIset{1}}};
			
			\node at (-0.1,0.9) {$ \Omega $}[grow cyclic absolute] 
			child[sibling angle=180,level distance=\lengthmap]{node{\DAWset{-1}}
				child[sibling angle=180,level distance=\lengthmap]{node{\DAWset{-2}}
					child[dashed,sibling angle=180,level distance=1.15*\lengthmapd]{node{}}
					child[dashed,sibling angle=180+atan(tan(\anglemapd)*1/4),level distance=1.15*\lengthmapd]{node{}}
					child[dashed,sibling angle=180+atan(tan(\anglemapd)*2/4),level distance=1.15*\lengthmapd]{node{}}
					child[dashed,sibling angle=180+atan(tan(\anglemapd)*3/4),level distance=1.15*\lengthmapd]{node{}}
					child[dashed,sibling angle=180+\anglemapd,level distance=1.15*\lengthmapd]{node{}}
				}
				child[sibling angle=180+atan(tan(\anglemap)*1/4),level distance=sqrt(\lengthmap^2+(\lengthmap*tan(\anglemap)*1/4)^2)]{node{\DAXset{-2}}}
				child[sibling angle=180+atan(tan(\anglemap)*2/4),level distance=sqrt(\lengthmap^2+(\lengthmap*tan(\anglemap)*2/4)^2)]{node{\DAKset{-2}}}
				child[sibling angle=180+atan(tan(\anglemap)*3/4),level distance=sqrt(\lengthmap^2+(\lengthmap*tan(\anglemap)*3/4)^2)]{node{\DAMset{-2}}}
				child[sibling angle=180+\anglemap,level distance=\lengthmap/cos(\anglemap)]{node{\DAIset{-2}}}
			}
			child[sibling angle=180+atan(tan(\anglemap)*1/4),level distance=sqrt(\lengthmap^2+(\lengthmap*tan(\anglemap)*1/4)^2)]{node{\DAXset{-1}}}
			child[sibling angle=180+atan(tan(\anglemap)*2/4),level distance=sqrt(\lengthmap^2+(\lengthmap*tan(\anglemap)*2/4)^2)]{node{\DAKset{-1}}}
			child[sibling angle=180+atan(tan(\anglemap)*3/4),level distance=sqrt(\lengthmap^2+(\lengthmap*tan(\anglemap)*3/4)^2)]{node{\DAMset{-1}}}
			child[sibling angle=180+\anglemap,level distance=\lengthmap/cos(\anglemap)]{node{\DAIset{-1}}};
			
			\draw[rounded corners=5pt,red,densely dotted] (0.91,0.2) rectangle (0.68,1);
			\node at (1.05,0.065) (nodeOmegam1B) {$ \hat{\Omega}_{1} $};
			\draw[red,densely dotted]  (1,0.065) -- (0.91,0.2);
			
			\draw[rounded corners=5pt,red,densely dotted] (1.61,0.17) -- (1.61,1) -- (1.38,1) -- (1.38,0.77) -- (0.65,0.77) -- (0.65,0.17) -- cycle;
			\node at (1.75,0.065) (nodeOmegam1B) {$ \hat{\Omega}_{2} $};
			\draw[red,densely dotted]  (1.7,0.065) -- (1.61,0.17);
			
			\draw[rounded corners=5pt,red,densely dotted] (-0.92,0.2) rectangle (-0.68,1);
			\node at (-1.05,0.065) (nodeOmegam1B) {$ \hat{\Omega}_{-1} $};
			\draw[red,densely dotted]  (-1.05,0.065) -- (-0.92,0.2);
			
			\draw[rounded corners=5pt,red,densely dotted] (-1.62,0.17) -- (-1.62,1) -- (-1.38,1) -- (-1.38,0.77) -- (-0.65,0.77) -- (-0.65,0.17) -- cycle;
			\node at (-1.78,0.065) (nodeOmegam1B) {$ \hat{\Omega}_{-2} $};
			\draw[red,densely dotted]  (-1.78,0.065) -- (-1.62,0.17);
			\end{tikzpicture}
		\end{scaletikzpicturetowidth}
	\end{center}
	\caption{Sets relations according to classification with DA and ADS.}
	\label{fig:graphDA}
\end{figure}
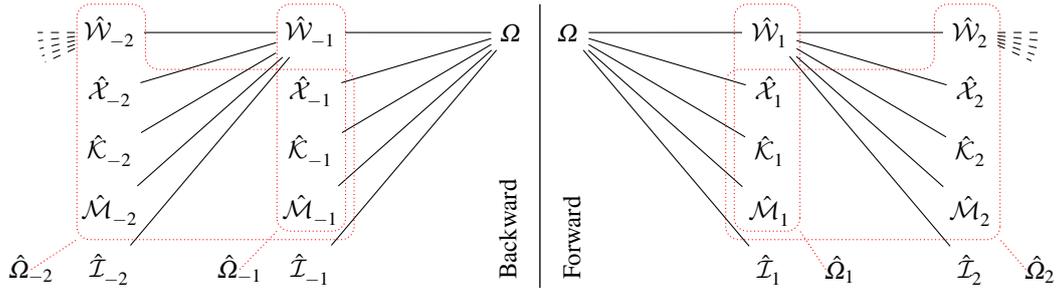

Therefore, the definitions of \Subsec{sec:set_defs} are adapted to \gls*{DA} propagation as follows:
\begin{enumerate*}[label=\roman*)]
    \item \textit{inconsistent} (sub-set \DAIset{i}) if the sub-domain performed the maximum number of splits allowed and tries to split again before completing the $i$-th period;
    \item \textit{weakly stable} (sub-set \DAWset{i}) if the sub-domain is consistent, and performs $i$ complete periods without escaping or impacting the target or its moons;
    \item \textit{unstable} (sub-set \DAXset{i}) if the sub-domain is consistent, and escapes from the target before completing the $i$-th;
    \item \textit{target--crash} (sub-set \DAKset{i}) if the sub-domain is consistent, and impacts with the target before completing the $i$-th period;    
    \item \textit{moon--crash} (sub-set \DAMset{i}) if the sub-domain is consistent, and impacts with one of the target's moons before completing the $i$-th period.
\end{enumerate*}
Conditions i), and ii)-v) apply after the particle performs $(i-1)$ periods around the target. As in \Subsec{sec:set_defs}, the sub-sets are defined for $ i \in \mathbb{Z} \textbackslash \{0\}  $ and the same considerations about propagation direction still apply. Moreover, the domain that is consistent after $i$ periods is defined as
\begin{equation}
    \hat{\Omega}_i = \Omega \textbackslash \bigcup_{j=1}^{\vert i \vert} \hat{\mathcal{I}}_{\sgn{i}j}.
	\label{eq:def_omega_set}
\end{equation}
In \Fig{fig:graphDA}, a graph reporting the relations between the sub-sets adapted to \gls*{DA} propagation is shown. As usual, a capture set is defined in \gls*{DA} propagation as $ \hat{\mathcal{C}}_{-1}^{n} := \hat{\mathcal{W}}_{n} \cap\hat{ \mathcal{X}}_{-1} $.

\subsection{Classification algorithm} \label{sec:classification}
The classification algorithm shown in \Alg{algo:classification_algorithm} is a while loop divided in two parts:
\begin{enumerate}
    \item the \gls*{DA} propagation of the stable set \DAWset{i - 1} from the current period to the next one, thus from $ T_{i-1} $ to $ T_{i} $;
    \item the extraction of all the sub-domains to classify them in the right set.
\end{enumerate}
A bridge is built between the definitions of \Subsec{sec:set_defs} and the ones of \Subsec{sec:hatted_sets}. The similarities are
\begin{equation}
	\hat{\Omega}_i \rightarrow \Omega, \quad
	\mathcal{\hat{I}}_{i} \rightarrow \emptyset, \quad
	\mathcal{\hat{K}}_{i} \rightarrow \mathcal{K}_{i}, \quad
	\mathcal{\hat{M}}_{i} \rightarrow \mathcal{M}_{i}, \quad
	\mathcal{\hat{X}}_{i} \rightarrow \mathcal{X}_{i}, \quad
	\mathcal{\hat{W}}_{i} \rightarrow \mathcal{W}_{i}, \quad
	\mathcal{\hat{C}}_{-1}^{i} \rightarrow \mathcal{C}_{-1}^{i}, \quad
	\emptyset \rightarrow \mathcal{D}_{i}.
	\label{eq:point-wise_vs_DA}    
\end{equation}
Note that the inconsistent set \DAIset{i} is mapped to the empty set due to the inability to compute these sub-domains with accuracy. Moreover, the acrobatic set \Dset{i} is not mapped by the \gls*{DA} classification algorithm. Nonetheless, this region of the search space represents a small fraction of the total for long enough propagation times \citep{luo2014constructing}.

\begin{algorithm}[tbp]
	\SetAlgoLined
	Set either $ i = 1 $ (forward propagation) or $ i = -1 $ (backward propagation)\;
	Set the search space \DAWset{0} $ = \Omega $\;
	Set capture epoch $ T_{0} = t_{0} $\;
	Set maximum number of periods $ n_{\mathrm{max}} $\;
	\While{$ |i| \leq n_{\mathrm{max}} $}{
		Propagate \DAWset{i - \sgn{i}} from $T_{i-\sgn{i}}$ to $T_{i}$\; 
		Extract inconsistent sub-domains in \DAIset{i}\;
		Extract crash sub-domains in \DAKset{i}\;
		Extract moon-crash sub-domains in \DAMset{i}\;
		Extract escaped sub-domains in \DAXset{i}\;
		Retain remaining sub-domains as \DAWset{i}\;
		Set $ i = i + \sgn{i} $\;
	}
	\caption{Classification algorithm.}
	\label{algo:classification_algorithm}
\end{algorithm}

\subsection{Propagation time span} \label{sec:timespan}
\begin{figure}[tbp]
\centering
\begin{subfigure}{.5\textwidth}
  \centering
  \includegraphics[width=\linewidth]{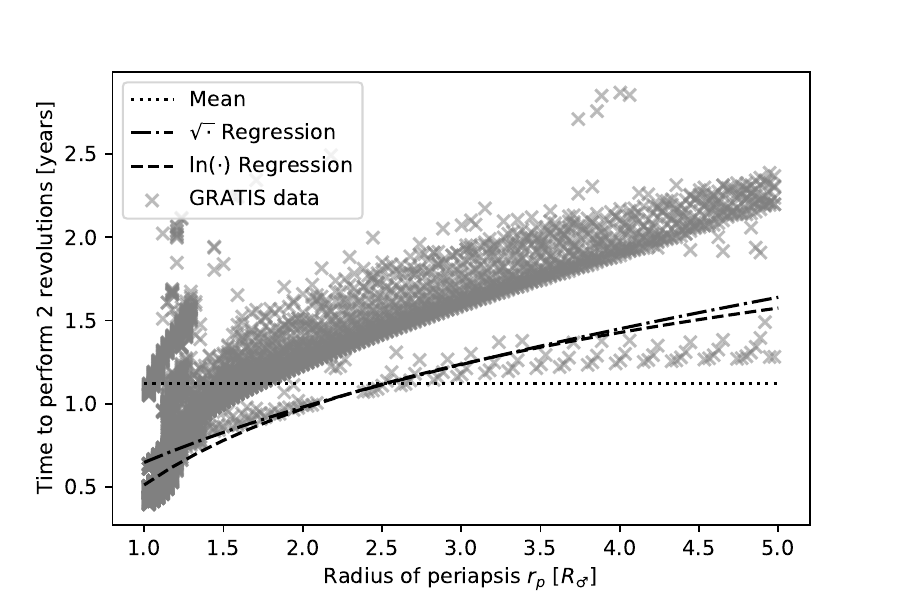}
  \caption{After 2 revolutions.}
\end{subfigure}%
\begin{subfigure}{.5\textwidth}
  \centering
  \includegraphics[width=\linewidth]{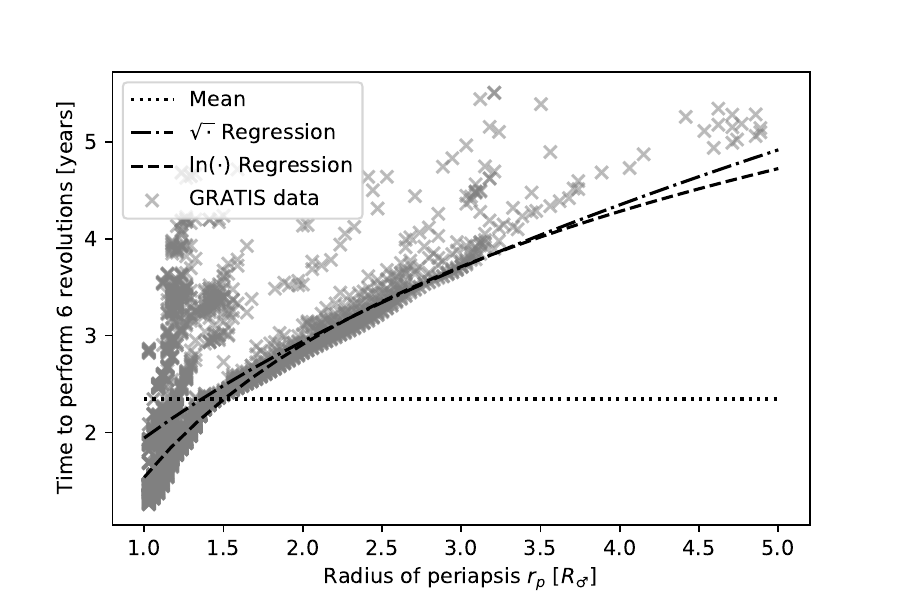}
  \caption{After 6 revolutions.}
\end{subfigure}
\caption{Regression of the revolution periods as a function of the radius of the periapsis $ r_{p} $. Data from point-wise computations issued by \gls*{GRATIS}.}
\label{fig:regression}
\end{figure}

Since it is not possible to track revolutions geometrically, they are tracked by counting the revolution periods. Thus, the need to estimate them with fidelity. These periods are determined by a least-square regression on data from point-wise computations issued by \gls*{GRATIS}\footnote{The time regressions have been derived from \glspl*{IC} propagated taking into account also the gravitational attractions of Uranus (B), and Neptune (B), later discarded due to their negligible influence.}, with respect to the radius of the periapsis $ r_{p} $. Two regression shapes are chosen:
\begin{enumerate}
    \item a square root shape
    \begin{equation}
        \label{eq:sqrt_period}
        T_{\sqrt{\cdot}} = A + B \sqrt{r_{p}}, \ \text{with } A = -0.19329, \ B = 2.10555;
    \end{equation}
    \item a logarithmic shape
    \begin{equation}
        \label{eq:ln_period}
        T_{\ln{(\cdot)}} = A + B \ln{(r_{p})}, \ \text{with } A = -7.35410, \ B = 1.44254.
    \end{equation}
\end{enumerate}
\Fig{fig:regression} shows these regressions for 2 and 6 periods, compared to the dataset generated by \gls*{GRATIS}.

\begin{figure}[tbp]
\centering
\begin{subfigure}{.5\textwidth}
  \centering
  \includegraphics[width=\linewidth]{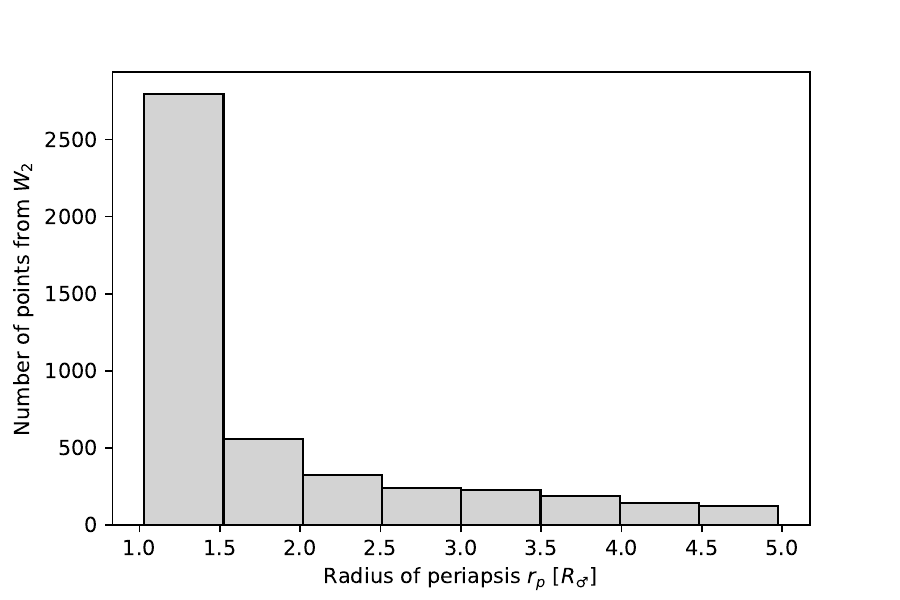}
  \caption{After 2 revolutions.}
\end{subfigure}%
\begin{subfigure}{.5\textwidth}
  \centering
  \includegraphics[width=\linewidth]{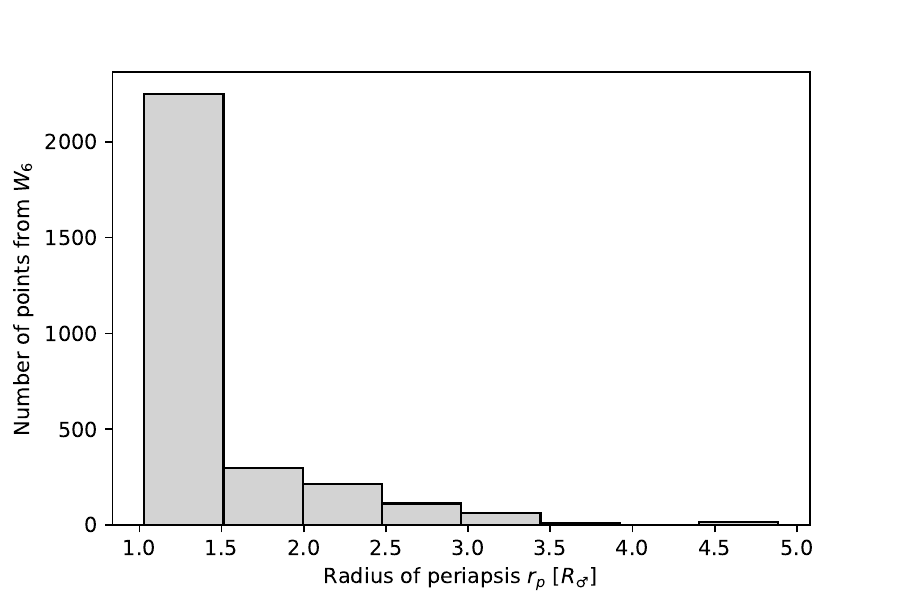}
  \caption{After 6 revolutions.}
\end{subfigure}
\caption{Histogram of the periapsis $ r_{p} $. Data from point-wise computations issued by \gls*{GRATIS}.}
\label{fig:histograms}
\end{figure}

\begin{figure}[tbp]
\centering
\begin{subfigure}{.5\textwidth}
  \centering
  \includegraphics[width=\linewidth]{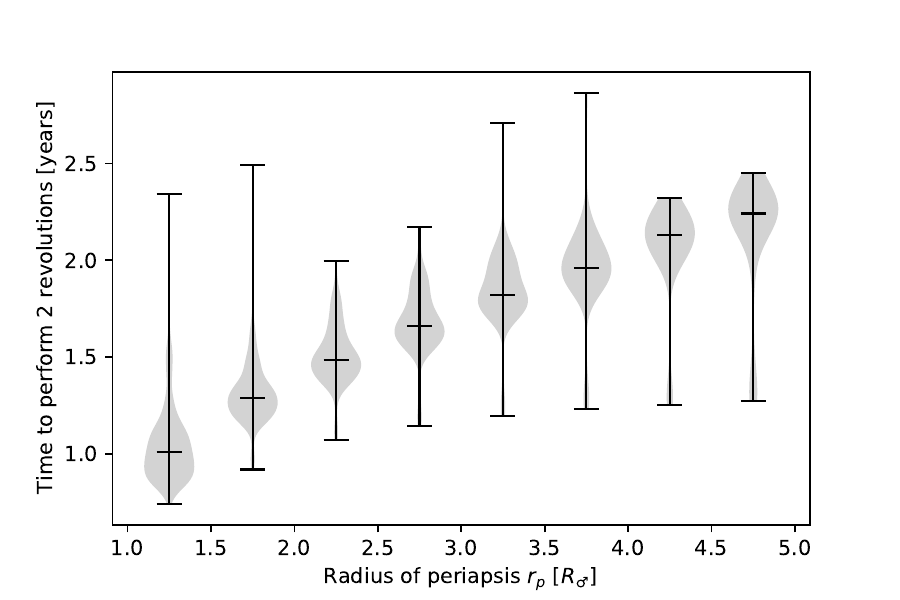}  \caption{After 2 revolutions.}
\end{subfigure}%
\begin{subfigure}{.5\textwidth}
  \centering
  \includegraphics[width=\linewidth]{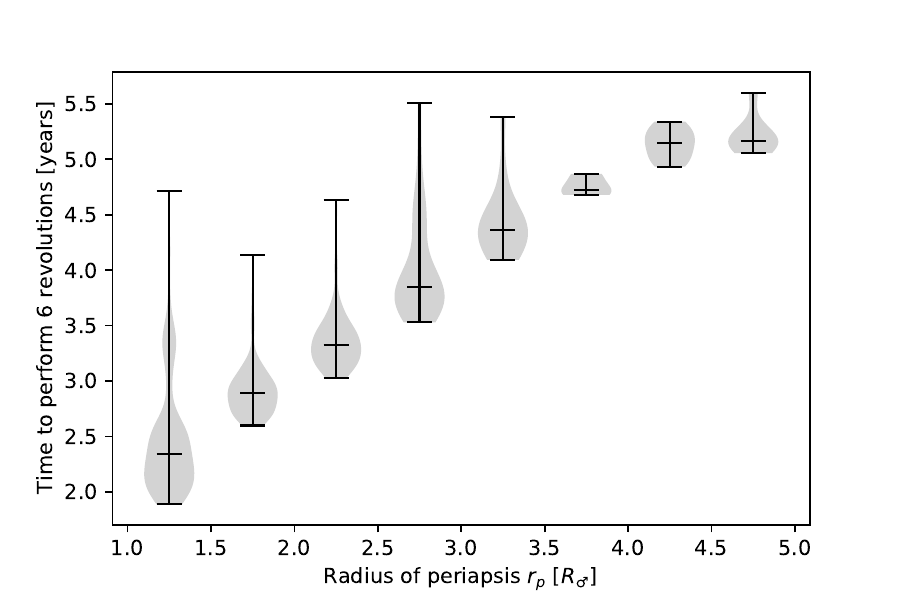}  \caption{After 6 revolutions.}
\end{subfigure}
\caption{Distribution of the revolution periods as a function of the radius of the periapsis $ r_{p} $. Data from point-wise computations issued by \gls*{GRATIS}.}
\label{fig:violin_plots}
\end{figure}

For the rest of this work, the logarithmic shape is chosen over the square root one, due to a better fit. Nevertheless, the regression does not fit well for 2 revolutions at large $r_p$, due to the high concentration of points at small radii, as shown by \Fig{fig:histograms}. Moreover, \Fig{fig:violin_plots} highlights how the distribution of the revolution times is widely spread, with respect to $r_p$, as the whiskers represent the minimum and maximum values, with the median in the middle. Therefore, the regressions displayed in \Eqs{eq:sqrt_period} and~\eqref{eq:ln_period} make up for a strong hypothesis on the behavior on the revolution periods.

\subsection{Consistency and quality criteria} \label{sec:criteria}
Once the classification of the search space is performed, tools need to be developed to evaluate it. Two criteria are defined to do so. The first is the \textit{consistency} criterion, which assesses the parts of the domain where the \gls*{DA} mapping cannot be trusted. The second is the \textit{quality} criterion, which assesses the performances of the \gls*{DA} mapping with respect to point-wise reference mapping carried out with \gls*{GRATIS}.

\paragraph{Consistency criterion.} \label{sec:consistency_def}
The consistency is a number assigned to each sub-domain. After the classification algorithm is performed, each sub-domain is given a consistency of either 1 if it does not belong to an inconsistent set \DAIset{i}, or 0 if it does. Consequently, the default consistency of a sub-domain is 1. Yet, if a sub-domain reaches the minimum size allowed and tries to split again, it is set to 0. The main interest is to compute the global consistency on all of the search space. It can be computed by performing the mean of the consistencies, pondered by the size of each sub-domain. The global consistency will then represent the ratio of the mapping that is not inconsistent. Therefore, it is the proportion of the mapping where \gls*{ADS} guarantees the accuracy of the mapping.

\paragraph{Quality criterion.} \label{sec:quality_def}
The quality criterion is a value assigned to each set issued from \gls*{GRATIS}'s classification. It represents the proportion of the set from \gls*{GRATIS} that is well-mapped by the \gls*{DA} classification algorithm. To be computed, the quality criterion requires a sample propagated point-wise with \gls*{GRATIS}. The quality $ q_\mathcal{A} $ of a \gls*{GRATIS} set $ \mathcal{A} $ is computed as follows:
\begin{enumerate}
    \item for each point $ x \in \mathcal{A} $, if $ x $  belongs to the set $ \mathcal{\hat{A}} $, according to the bridge between \gls*{DA} and point-wise classification exposed in \Eq{eq:point-wise_vs_DA}, then the quality of that point $ q_{x} = 1 $, otherwise $ q_{x} = 0 $;
    \item the quality $q_\mathcal{A}$ of a set $\mathcal{A}$ is the mean of all the values of $\left(q_x\right)_{x\in \mathcal{A}}$. In other words
    \begin{equation}
        q_\mathcal{A} = \frac{\#\left(\mathcal{A}\cap\mathcal{\hat{A}}\right)}{\#\left(\mathcal{A}\right)}.
    \end{equation}
    In the case where $\#\left(\mathcal{A}\right)=0$, then $\#\left(\mathcal{A}\cap\mathcal{\hat{A}}\right)=0$. Thus, by convention $q_\mathcal{A} = 1$;
    \item the confidence intervals on $ q_\mathcal{A} $ are evaluated according to \citet{robert2004monte}, and \citet{hanley1983if}.
\end{enumerate}
The quality represents the probability for a point from the point-wise mapping to be mapped correctly with the \gls*{DA} mapping. Note that since \DAIset{i} is mapped to $\emptyset$, the quality criterion of $\Omega_i$ is less than or equal to the consistency criterion after $i$ revolutions. A schematic representation clarifying the meaning of the quality criterion is shown in \Fig{fig:quality_criterion_scheme}.

\begin{figure}
    \centering
    \begin{tikzpicture}
    
    \filldraw[fill=blue!15!white, draw=black] (0,0) rectangle (4,4);
    \draw (1,0) arc (0:90:1cm) node[right=0.4cm, below=0.3cm]{$\hat{\mathcal{A}}$};
    \draw (1.2,1.3) node[blue] {$\times$} node[above=0.1cm] {$\mathcal{A}$}; %
    \draw (3.7,2) node[green] {$\circ$} node[above=0.1cm] {$\mathcal{B}$}; %
    
    \filldraw[fill=green!15!white, draw=black] (0,4) rectangle (4,8);
    \draw (1,4) arc (0:90:1cm) node[right=0.4cm, below=0.3cm]{$\hat{\mathcal{B}}$};
    \draw (0.5,5.1) node[green] {$\times$} node[above=0.1cm] {$\mathcal{B}$}; %
    \draw (2.5,7) node[green] {$\times$} node[above=0.1cm] {$\mathcal{B}$}; %
    
    \filldraw[fill=red!15!white, draw=black] (4,0) rectangle (8,8);
    \draw (5,0) arc (0:90:1cm) node[right=0.4cm, below=0.3cm]{$\hat{\mathcal{C}}$};
    \draw (5.2,1.1) node[blue] {$\circ$} node[above=0.1cm] {$\mathcal{A}$}; %
    \draw (7.7,4) node[red] {$\times$} node[above=0.1cm] {$\mathcal{C}$}; %
    \draw (4.3,6.7) node[green] {$\circ$} node[above=0.1cm] {$\mathcal{B}$}; %
    \draw (5.2,4.8) node[red] {$\times$} node[above=0.1cm] {$\mathcal{C}$}; %
    
    \draw (0,-1) node {Legend:};
    
    \filldraw[fill=black!15!white, draw=black] (0,-3) rectangle (1,-2);
    \draw (1,-2.5) node[right] {Sub-domain classified by \gls*{DA}};
    
    \draw (7.5,-2) node[black] {$\times$} node[right=0.5cm] {\gls*{GRATIS} IC accurately mapped by \gls*{DA}}; %
    \draw (7.5,-3) node[black] {$\circ$} node[right=0.5cm] {\gls*{GRATIS} IC inaccurately mapped by \gls*{DA}}; %
    
    \draw (11,7) node {$q_\mathcal{A} = \frac{1}{2}=0.5$}; %
    \draw (11,5) node {$q_\mathcal{B} = \frac{2}{4}=0.5$}; %
    \draw (11,3) node {$q_\mathcal{C} = \frac{2}{2}=1$}; %
    \draw (11,1) node {$\#\left(\mathcal{D}\right)=0 \implies q_\mathcal{D} = 1$}; %
    
    \end{tikzpicture}
    \caption{Example of computation of the quality criterion for a given mapping.}
    \label{fig:quality_criterion_scheme}
\end{figure}

\subsection{Representation methods} \label{sec:representation}
The representation of the results produced by this methodology raises two issues. Firstly, due to \gls*{ADS}, the largest fully-split sub-domains have a size of $ \approx \SI{50}{\kilo\meter}$, while the typical length of the search space is $\approx \SI{3000}{\kilo\meter} $. Thus, it is impossible to see them or smaller sub-domains on a global representation. To solve this problem, the resolution of the display is downgraded to a smaller one. Furthermore, instead of plotting each sub-domain independently, the density of sub-domains per pixel is preferred. It allows detecting areas with a high number of sub-domains even on a global visualization.

Secondly, the Cartesian representation of the search space tends to shrink details for low values of the radius of the periapsis $ r_{p} $, while the size of structures located at larges values of $ r_{p} $ is amplified. However, areas with a larger amount of data to visualize are located at small $ r_{p} $. Thus, results are visualized on the $ r_{p} \times \omega $ plane, which is the search space in Keplerian coordinates.

\section{Results} \label{sec:results}

In this work, the results from two simulations with initial grids having different resolution are presented. The order of the DA mappings is 20 to provide good accuracy on consistent domains and delay the triggering of \gls*{ADS}. Therefore, increasing the order lowered the computational time and enhanced consistency. Nevertheless, a larger order would not increase the approximation performances, see \citet{wittig2015propagation}. The first simulation relies on a coarse grid providing a low-resolution mapping where the search space is initially divided into $ 32 \times 32 $ domains and the maximum number of splits allowed to the \gls*{ADS} algorithm is 9. On the contrary, the second is computed on a finer grid returning a high-resolution mapping. In the latter, the search space is divided into $ 128 \times 128 $ domains and 10 maximum splits are allowed. Results are visualized in the low-resolution mapping in Cartesian coordinates. Sub-domains densities and last step epochs are computed on the whole search space. Then, the \gls*{DA} classification outcome is rendered for qualitative analysis. Finally, a quantitative analysis of those mappings exploiting consistency and quality criteria follows.

\subsection{Mapping} \label{sec:mapping}
\begin{figure}[tbp]
	\centering
	\begin{subfigure}{.5\textwidth}
		\centering
		\includegraphics[width=1.2\linewidth]{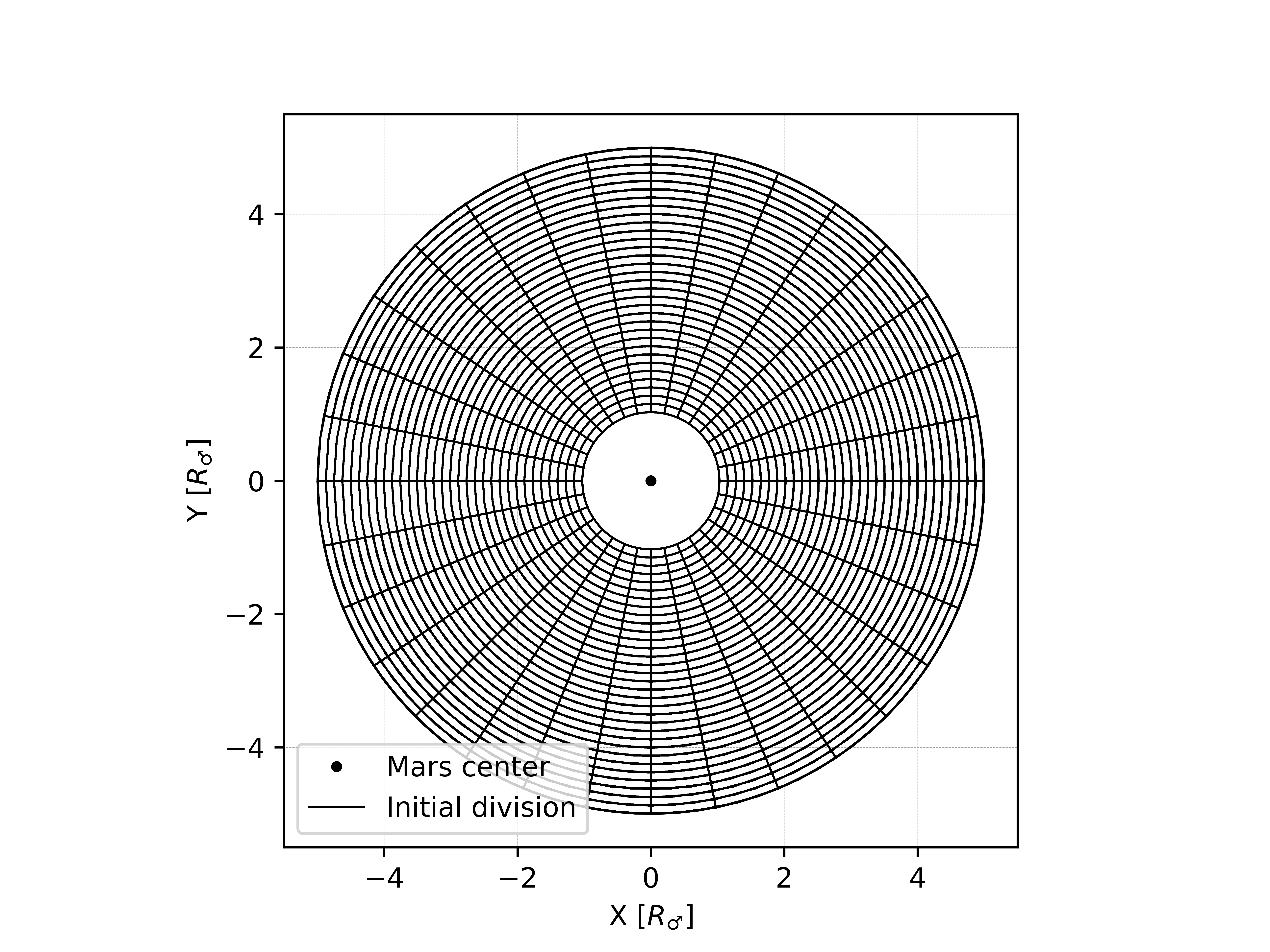}
		\caption{Before propagation with \gls*{ADS}.}
		\label{fig:car_map_left}
	\end{subfigure}%
	\begin{subfigure}{.5\textwidth}
		\centering
		\includegraphics[width=1.2\linewidth]{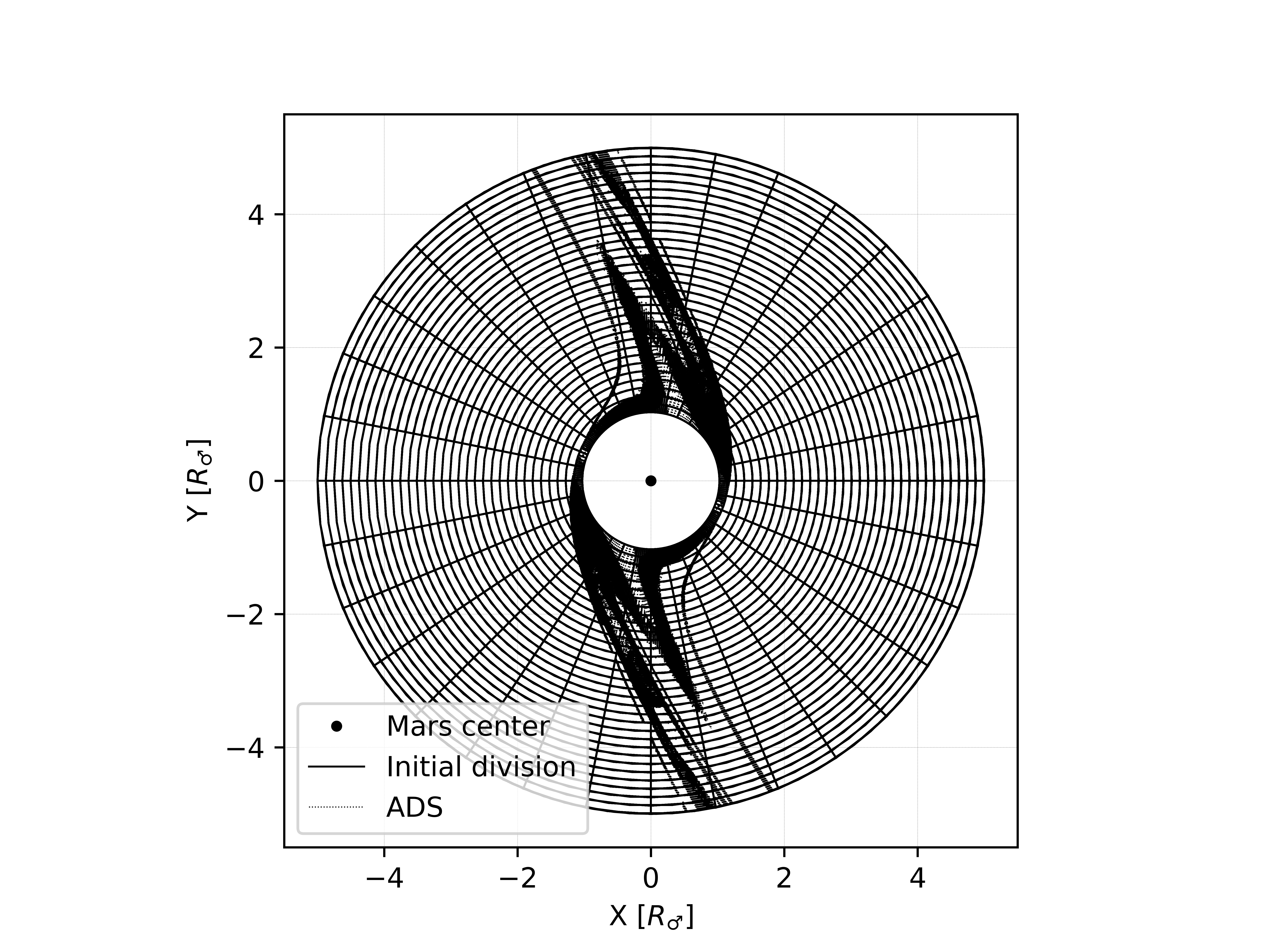}
		\caption{After propagation with \gls*{ADS}.}
		\label{fig:car_map_right}
	\end{subfigure}
	\caption{Cartesian visualization of the mapping of the search space.}
	\label{fig:cartesian_mapping}
\end{figure}

\Fig{fig:cartesian_mapping} represents the low-resolution search space in Cartesian coordinates, before and after the propagation with \gls*{ADS} (left and right, respectively). Before propagation, the search space is divided regularly into small domains for parallelization (see \Fig{fig:car_map_left}). Conversely, an irregular sub-divisions highlighting dynamical changes is shown in \Fig{fig:car_map_right}. This visualization shows the necessity to use different visualization methods. Indeed, it is hard to analyze close to zones where \gls*{ADS} creates a large number of domains, although these are the most interesting regions.

\subsubsection{Density of sub-domains and last step epochs} \label{sec:subdomain}
\begin{figure}[tbp]
	\centering
	\begin{subfigure}{.5\textwidth}
		\centering
		\includegraphics[width=1\linewidth]{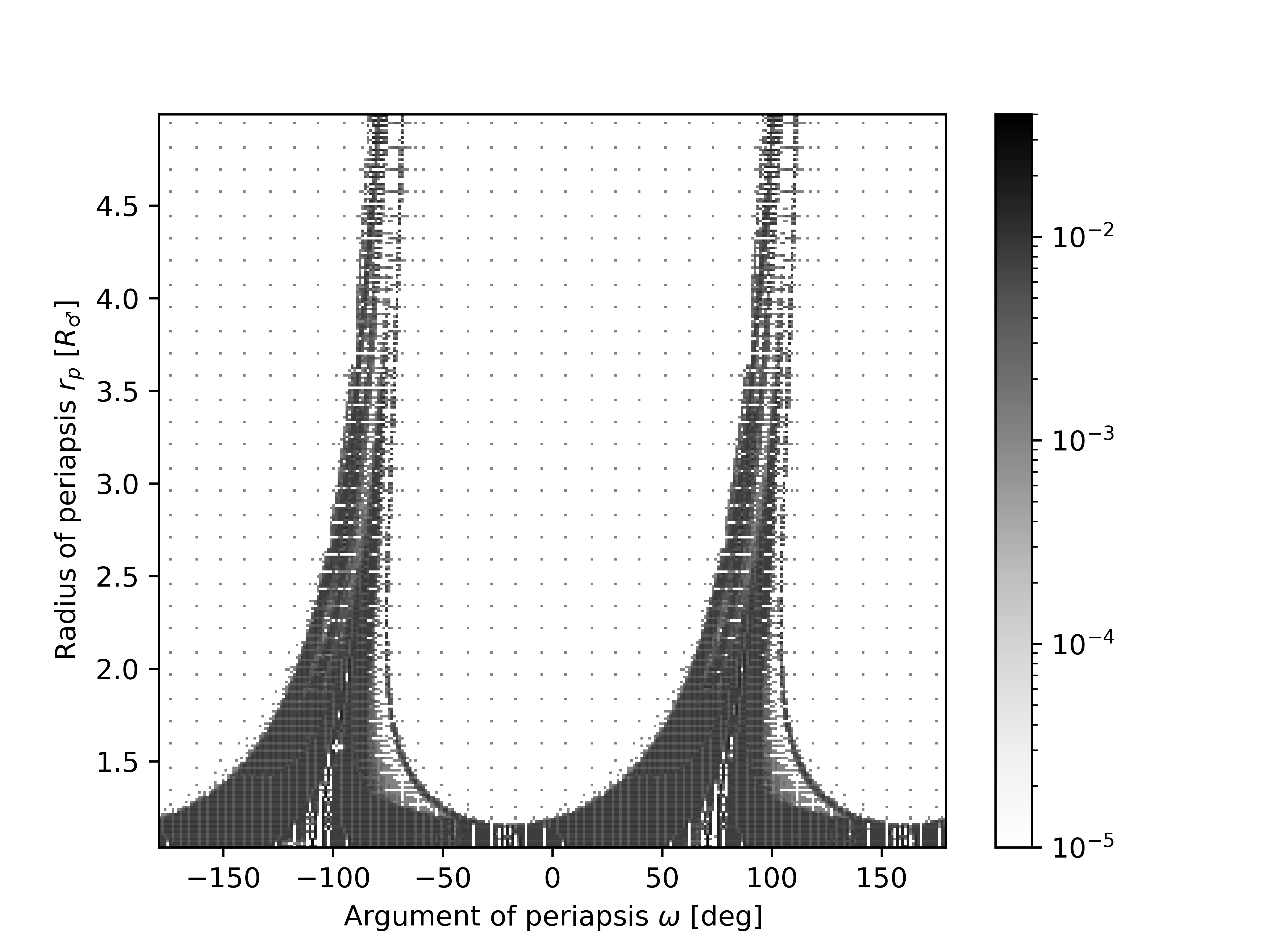}
		\caption{Low resolution.}
		\label{fig:density_mapping_left}
	\end{subfigure}%
	\begin{subfigure}{.5\textwidth}
		\centering
		\includegraphics[width=1\linewidth]{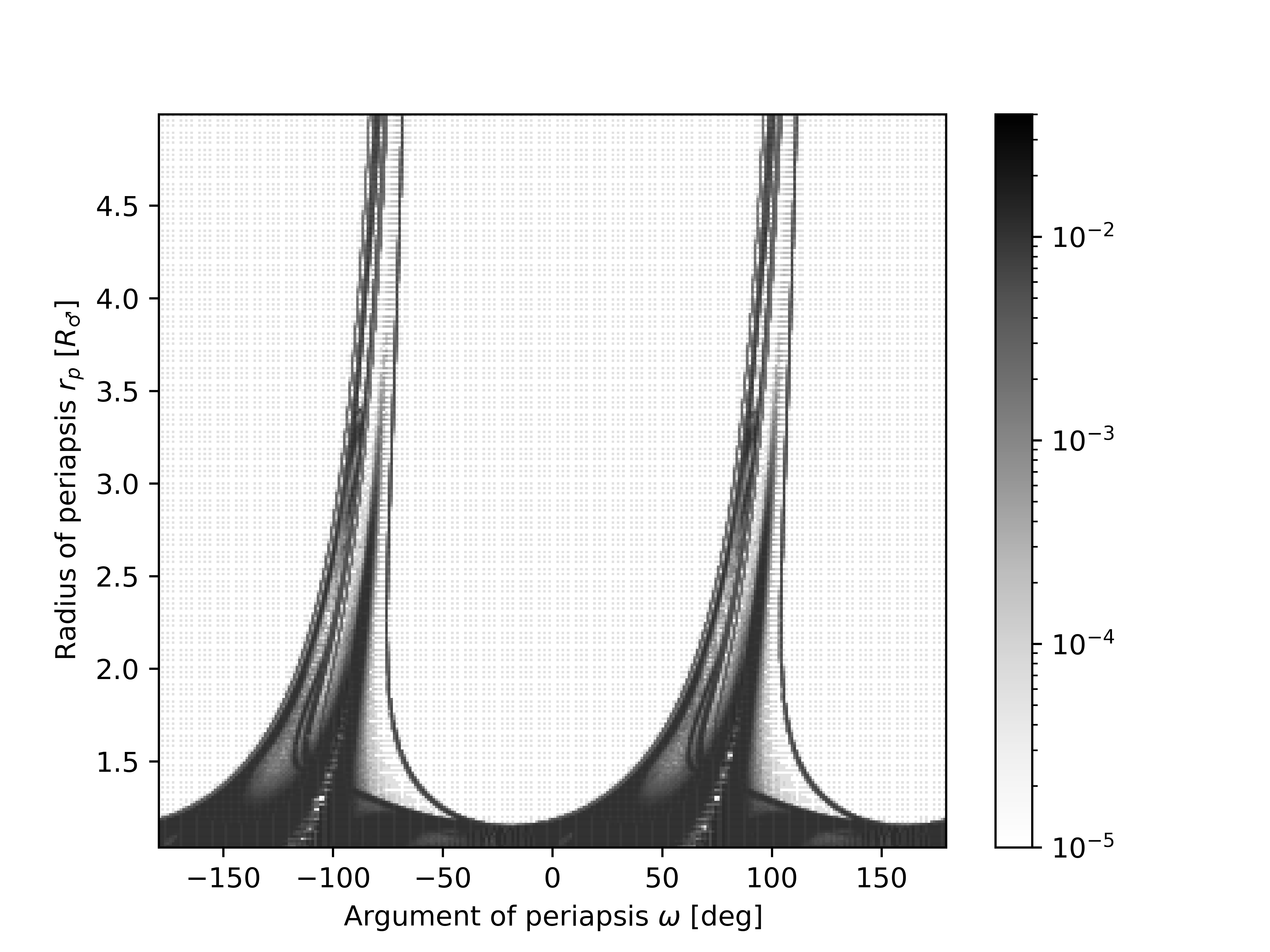}
		\caption{High resolution.}
		\label{fig:density_mapping_right}
	\end{subfigure}
	\caption{Density of sub-domains, in percent of the overall sum of sub-domains after 6 periods.}
	\label{fig:density_mapping}
\end{figure}

\Fig{fig:density_mapping} shows the density of sub-domains for the low-resolution and high-resolution mappings (left and right, respectively). They mostly differ because the high-resolution mapping provides more variations. Indeed, the density is either maximal (in grey) or minimal (in white) in the low-resolution one. On the contrary, the high-resolution mapping provides more shades, highlighting several degrees of nonlinearity accurately mapped. In addition, the low-resolution mapping created $7.2\times10^4$ sub-domains. It is fewer propagations than the point-wise cartography of \Sec{sec:classification_results} for the representation of the same search space. Meanwhile, the number of sub-domains in the high-resolution mapping is $1.2\times10^6$, thus, more than ten times the number of point-wise propagations. Therefore, there is no computational benefit in using high-resolution mapping instead of point-wise propagations.
However, as opposed to a point-wise mapping, the sub-domains are heterogeneously distributed, and are mainly located in highly nonlinear regions. For instance, the vicinity of Mars and areas that lead to collisions have a large sub-domains density compared to the majority of the search space. Furthermore, the sub-domain density draws the boundaries between different sets, due to their different dynamical behaviors.
On that account, \gls*{DA} mapping provides useful macroscopic information on the search space, while grid point-wise mapping does not allow highlighting highly nonlinear areas or boundaries.

\begin{figure}[tbp]
	\centering
	\begin{subfigure}{.5\textwidth}
		\centering
		\includegraphics[width=1\linewidth]{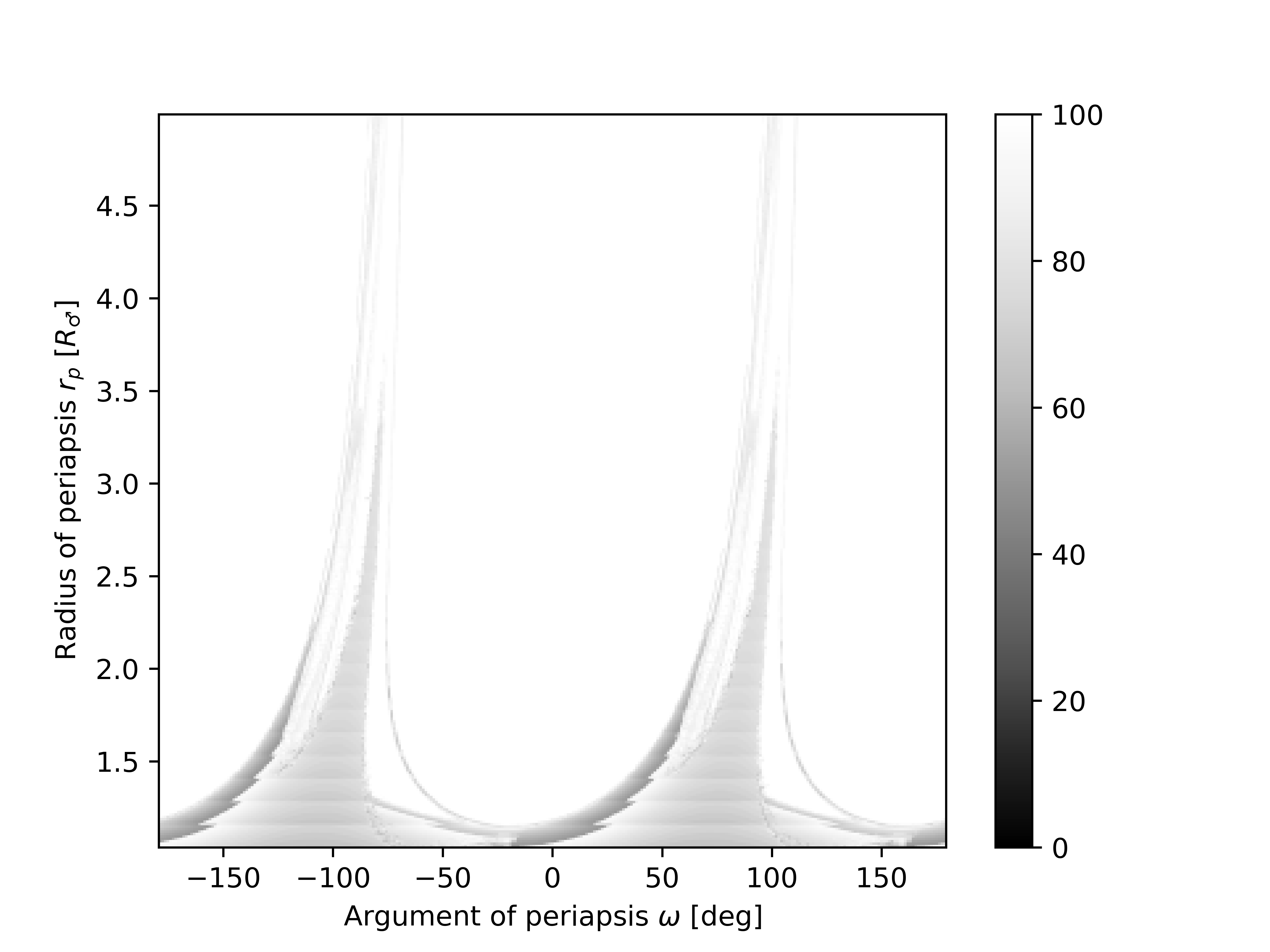} 
		\caption{Low resolution.}
		\label{fig:last_step_mapping_left}
	\end{subfigure}%
	\begin{subfigure}{.5\textwidth}
		\centering
		\includegraphics[width=1\linewidth]{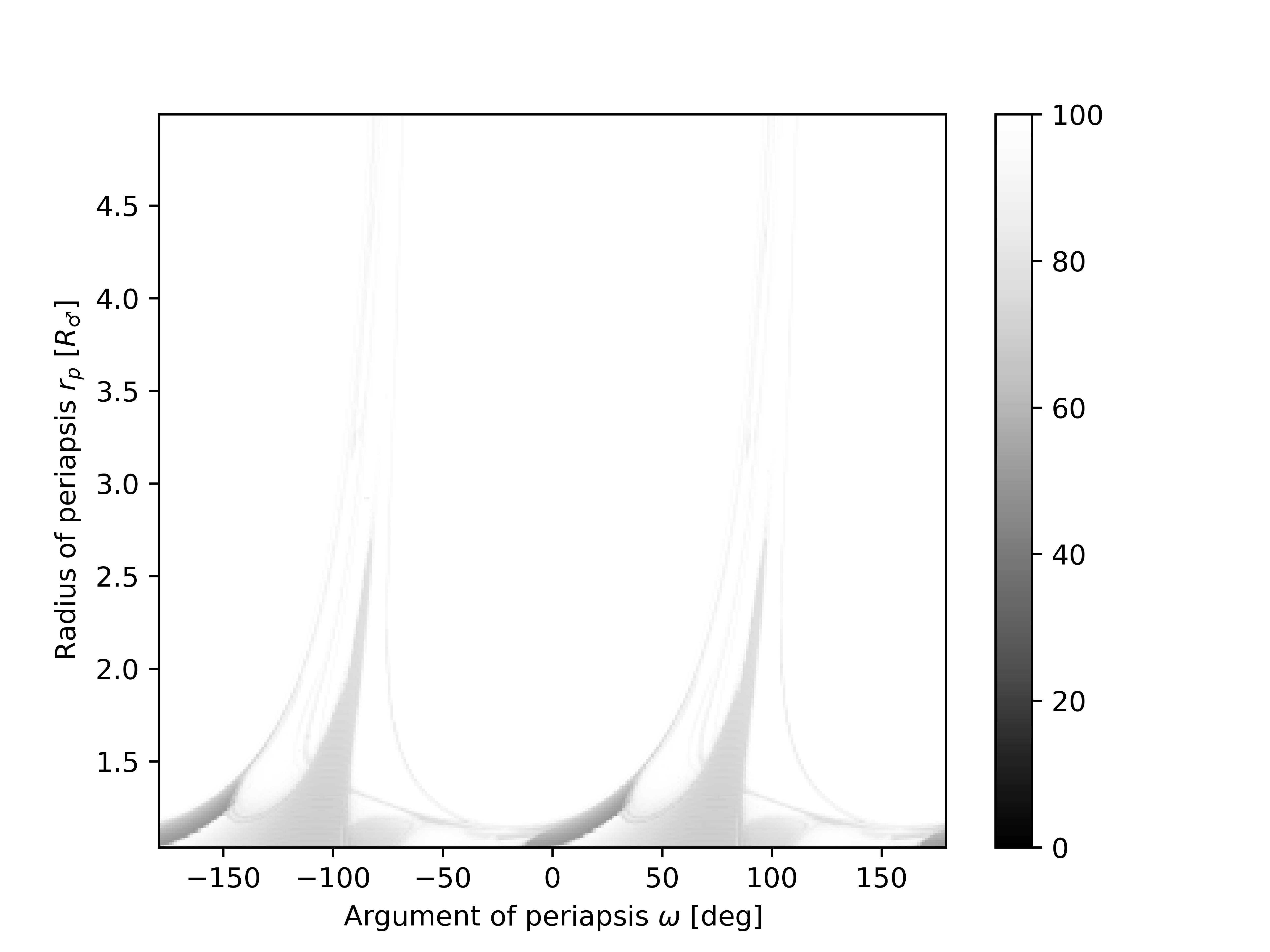}
		\caption{High resolution.}
		\label{fig:last_step_mapping_right}
	\end{subfigure}
	\caption{Last step epochs, in percent of the overall propagation time span corresponding to 6 periods.}
	\label{fig:last_step_mapping}
\end{figure}

Moreover, epochs of last steps carried out by the propagation scheme are shown in \Fig{fig:last_step_mapping}, in percent of the overall propagation time span. The integration of a sub-domain stops either after a collision occurs or when the sub-domain itself is declared inconsistent. Note that at low $ r_{p} $, the lightest zones on the low-resolution mapping (\Fig{fig:last_step_mapping_left}) become completely white when visualized in high resolution (\Fig{fig:last_step_mapping_right}). Therefore, propagation of these areas is now completed and sub-domains becomes consistent on a finer grid. In these regions, the consistency criterion is expected to rise when the mapping resolution increases.

\subsubsection{Classification} \label{sec:classification_results}
\begin{figure}[tbp]
	\centering
	\begin{subfigure}{.33\textwidth}
		\centering
		\includegraphics[width=1.1\linewidth]{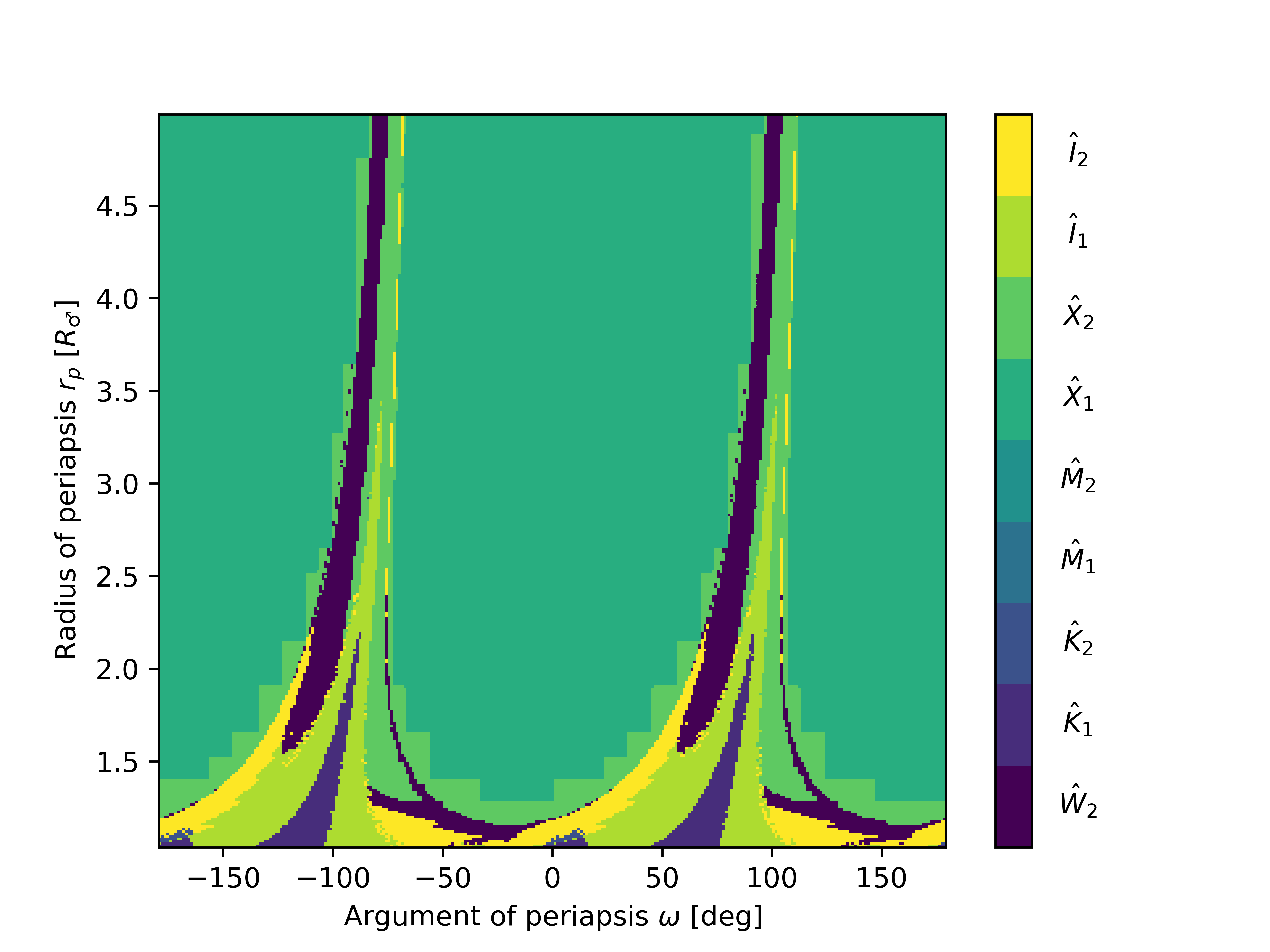}
		\caption{Low resolution.}
	\end{subfigure}%
	\begin{subfigure}{.33\textwidth}
		\centering
		\includegraphics[width=1.1\linewidth]{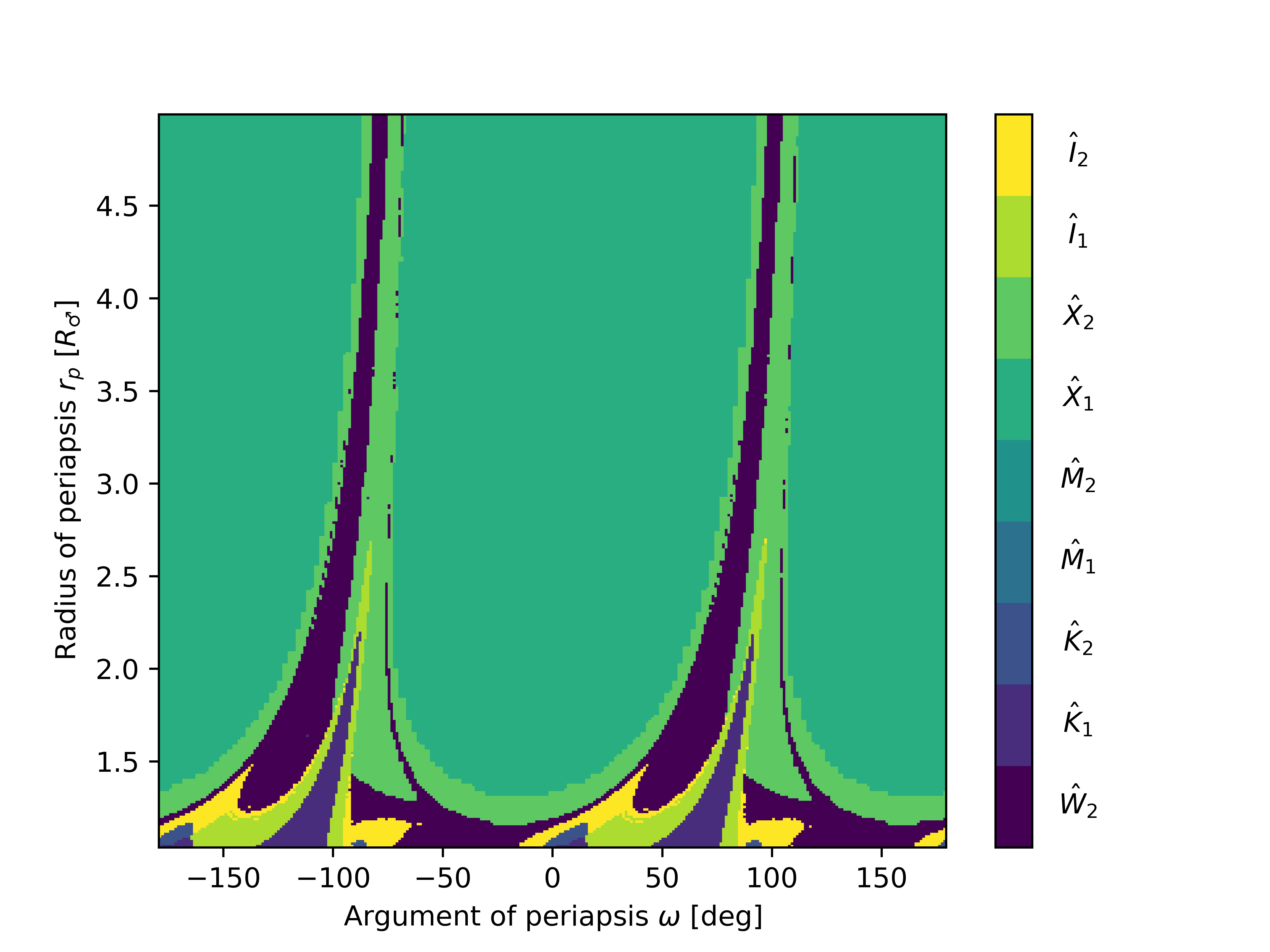}
		\caption{High resolution.}
	\end{subfigure}
	\begin{subfigure}{.33\textwidth}
		\centering
		\includegraphics[width=1.1\linewidth]{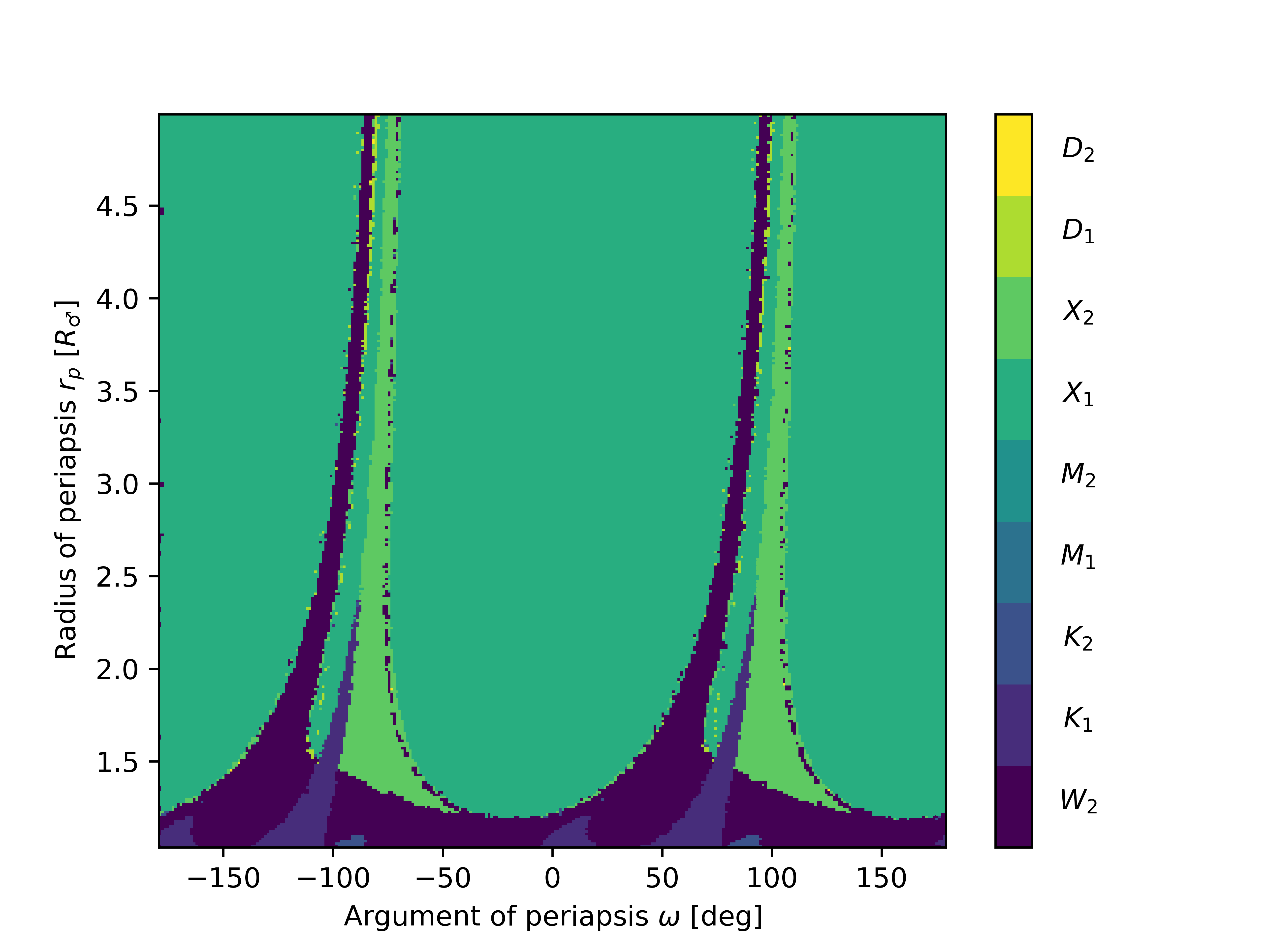}
		\caption{GRATIS.}
	\end{subfigure}
	\caption{Results of the classification algorithm after 2 revolutions compared to GRATIS.}
	\label{fig:classification_2_revolutions}
\end{figure}

\Fig{fig:classification_2_revolutions} shows the results of the classification algorithm after two revolutions for the low-resolution mapping (left), the high-resolution mapping (middle), and from point-wise propagation using \gls*{GRATIS} (right). The latter used as a reference and derived through point-wise propagation from a sample of $ 10^{5} $ points. Each color represents a different set, following the bridge between sets from \gls*{DA} to point-wise mapping established in \Sec{sec:classification}. Thus, a set \WSBset{A}{}{} from \gls*{GRATIS} classification is colored the same way as \WSBset{\hat{A}}{}{} from \gls*{DA} classification. Since inconsistent sets \DAIset{i} and acrobatic sets \Dset{i} have not corresponding sets in point-wise and DA mappings, respectively, they share the same colors.

While these three mappings have a similar macroscopic look, they differ in many ways. First of all, the sizes of the inconsistent sets decrease dramatically from the low-resolution mapping to the high-resolution one. Moreover, these gains in consistency mainly benefit the representation of the stable set \Wset{2}. Furthermore, the escape sets \DAXset{1} and \DAXset{2} are poorly mapped compared to the results of \gls*{GRATIS} on \Xset{1} and \Xset{2}. These issues in the mapping are mostly due to errors in the approximation of the revolution period.

\begin{figure}[tbp]
	\centering
	\begin{subfigure}{.33\textwidth}
		\centering
		\includegraphics[width=1.1\linewidth]{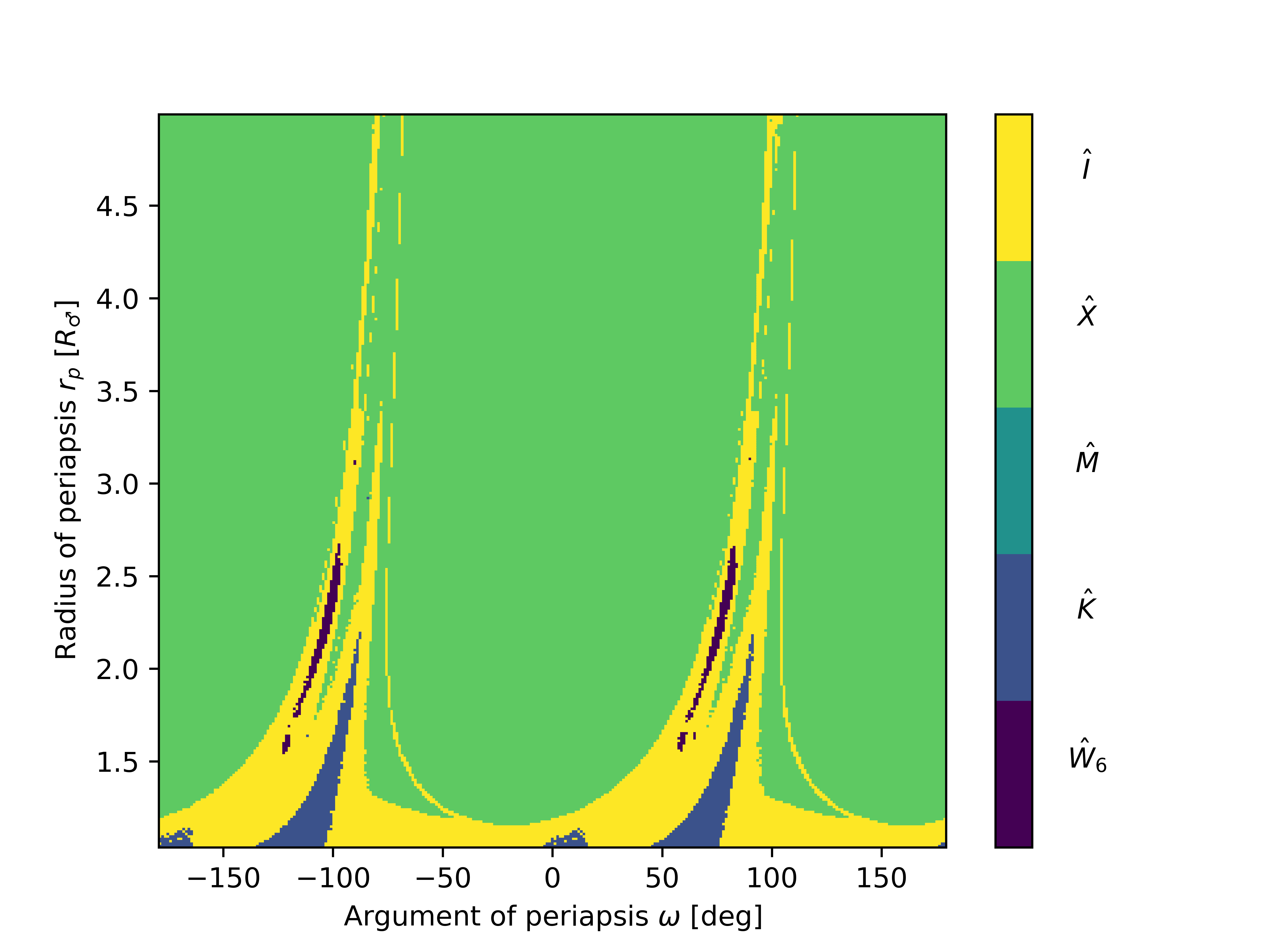}
		\caption{Low resolution.}
	\end{subfigure}%
	\begin{subfigure}{.33\textwidth}
		\centering
		\includegraphics[width=1.1\linewidth]{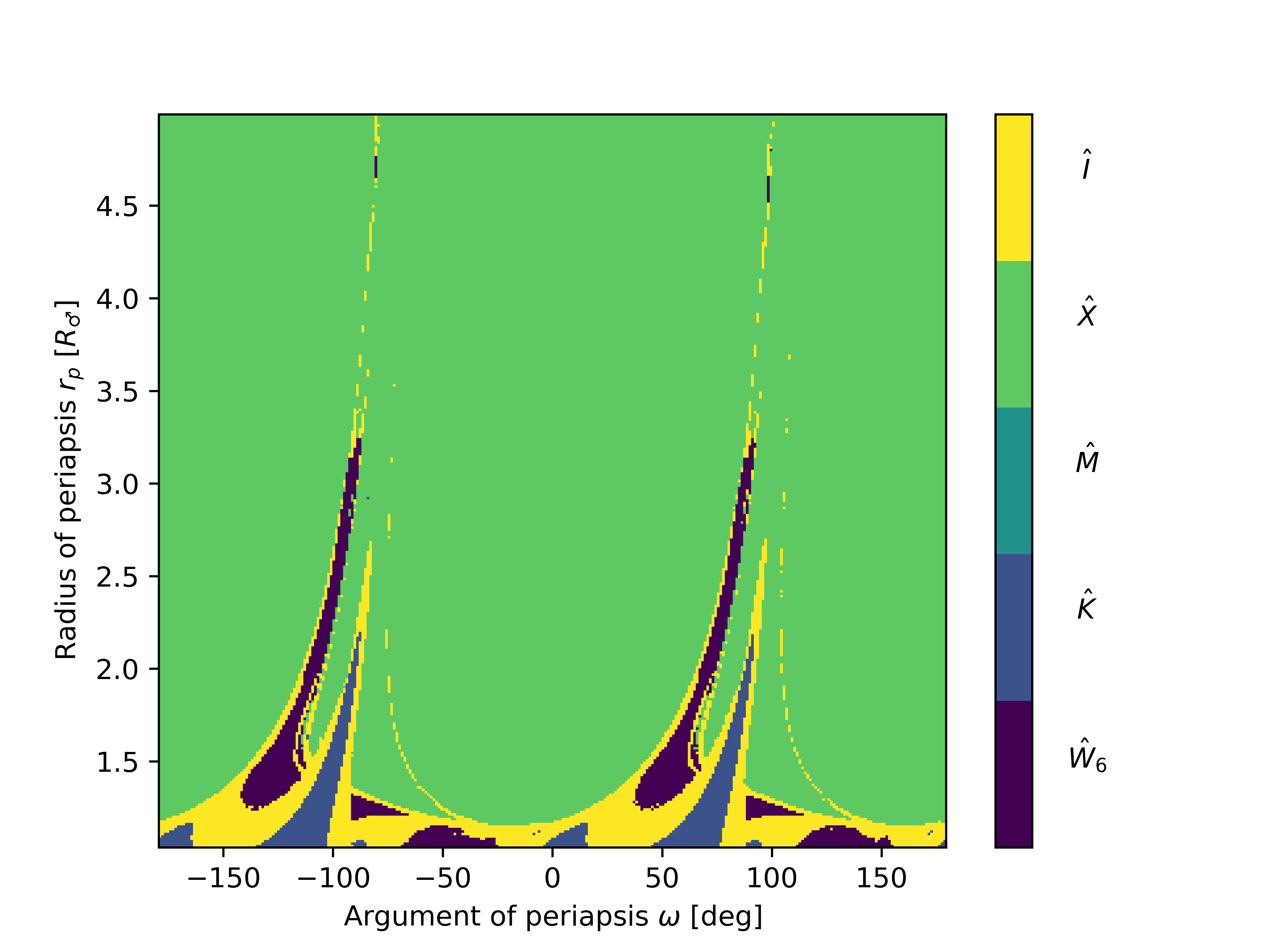}
		\caption{High resolution.}
	\end{subfigure}
	\begin{subfigure}{.33\textwidth}
		\centering
		\includegraphics[width=1.1\linewidth]{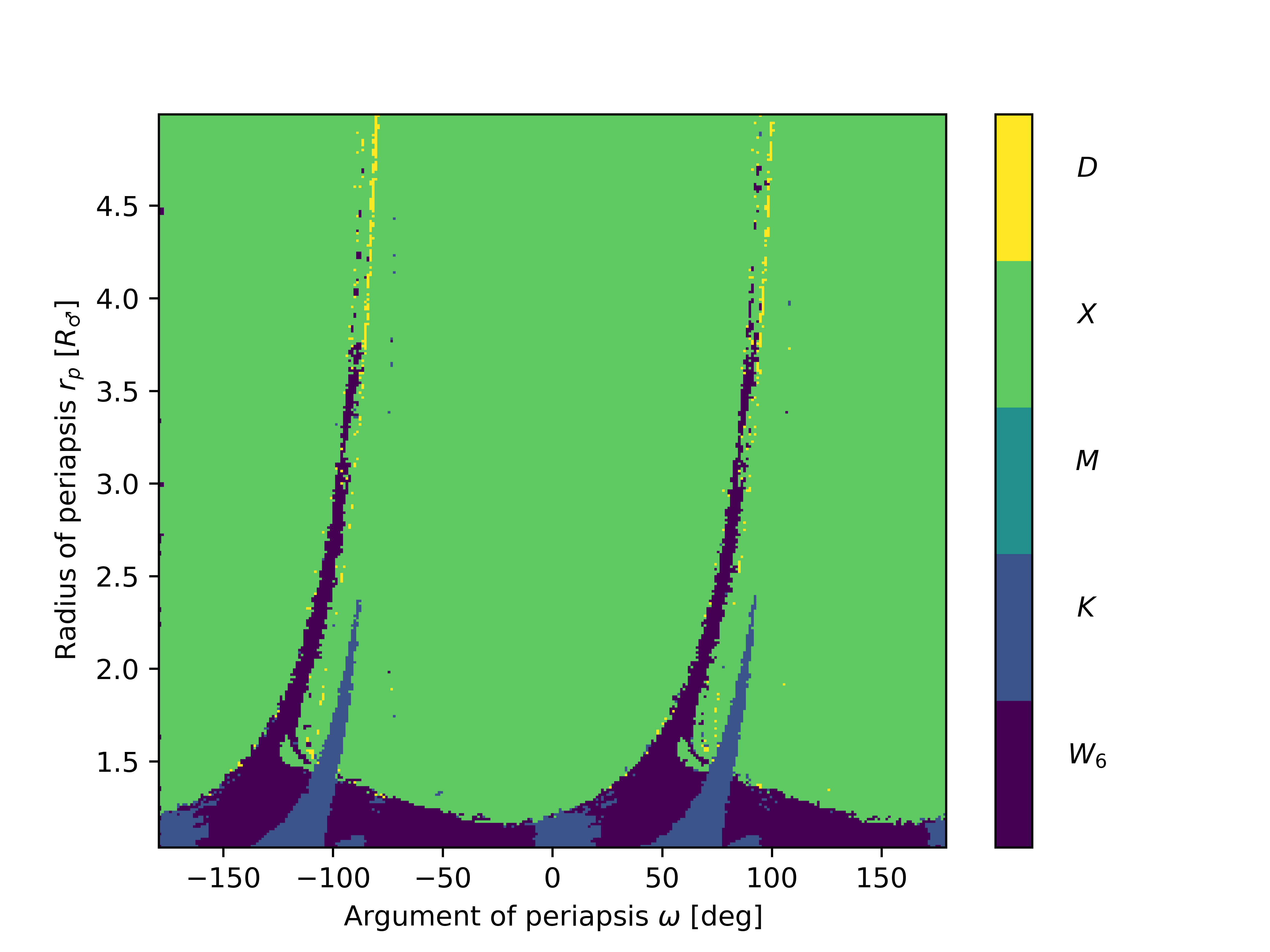}
		\caption{GRATIS.}
	\end{subfigure}
	\caption{Results of the classification algorithm after 6 revolutions compared to GRATIS.}
	\label{fig:classification_6_revolutions}
\end{figure}

\Fig{fig:classification_6_revolutions} presents the results of the classification algorithm after six revolutions for the low-resolution mapping (left), the high-resolution mapping (middle), and from point-wise propagation using \gls*{GRATIS} (right). Also in this case the latter mapping is used as reference and it is derived through point-wise propagation from a sample of $ 10^{5} $ points as before. These three charts show simplified mappings. For the sake of clarity, all sets $ \left( \mathcal{A}_{j} \right)_{j \in \llbracket 1,6\rrbracket}$ are united under the name $ \mathcal{A} $, defined as 
\begin{equation}
	\mathcal{A} = \bigcup_{j=1}^{6} {\mathcal{A}_{j}}.
	\label{eq:siplified_mapping_definition}
\end{equation}
This representation method is adopted for all sets but \Wset{6}, and \DAWset{6}.

The same observations can be done for Figs.~\ref{fig:classification_2_revolutions} and~\ref{fig:classification_6_revolutions}. However, \Fig{fig:classification_6_revolutions} shows a strong predominance of the inconsistent set, especially at low resolution. The high-resolution mapping improves these results but still struggles to handle long accurate propagations on such a large search space. Nevertheless, these results show that \gls*{DA} mapping with \gls*{ADS} can highlight changes of behavior all over the search space. Moreover, DA mapping provides a continuous description of the domain. The dynamical behavior is the same everywhere, providing consistency. Meanwhile, point-wise computations deliver discrete information challenging to interpolate due to nonlinearities.

\subsection{Quantitative performance analysis} \label{sec:results_criteria}
After performing a qualitative analysis on the resulting \gls*{DA} mappings, the two criteria developed deliver quantitative data on the computed mappings.

\subsubsection{Consistency criterion} \label{sec:consistency_results}
\begin{figure}[tbp]
	\centering
	\includegraphics[width=0.7\linewidth]{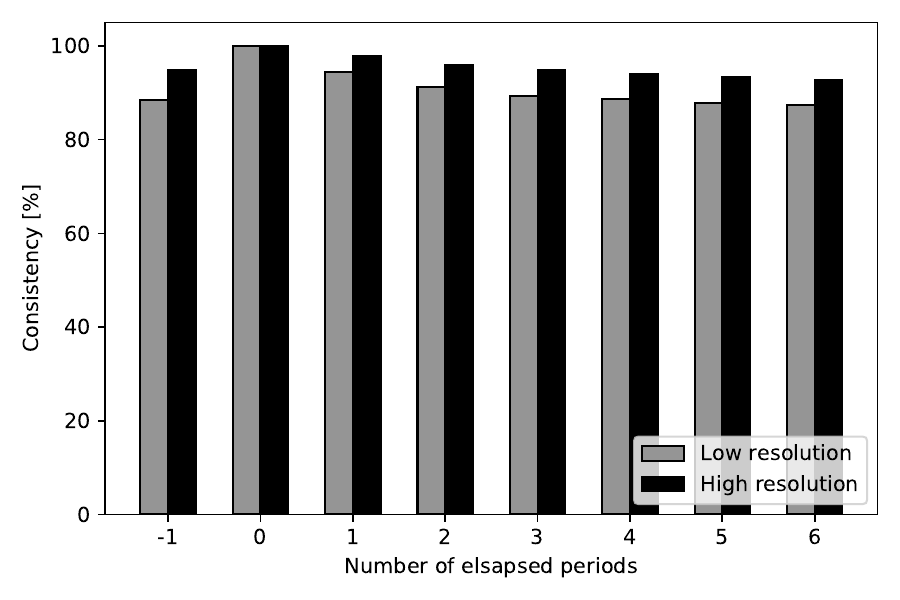}
	\caption{Consistency criterion for several number of revolutions and the two mapping resolutions.}
	\label{fig:consistency}
\end{figure}

In \Fig{fig:consistency}, the global consistency rate for both the low- and high-resolution mappings after several revolution periods are shown in bright and dark colors, respectively. The consistency at zero revolution is 100\% since mappings only represent \glspl*{IC}. It appears that the consistency remains high, even after six revolutions, since it stays above 87\%. It means that the size of the inconsistent sets is 13\% or less on the total search space. Moreover, the consistency of the high-resolution mapping is around 5 points above the one of the low-resolution one. It demonstrates the gain of accuracy delivered by these additional computations.

\subsubsection{Quality criterion} \label{sec:quality_results}
\begin{figure}[tbp]
	\centering
	\begin{subfigure}{.5\textwidth}
		\centering
		\includegraphics[width=1\linewidth]{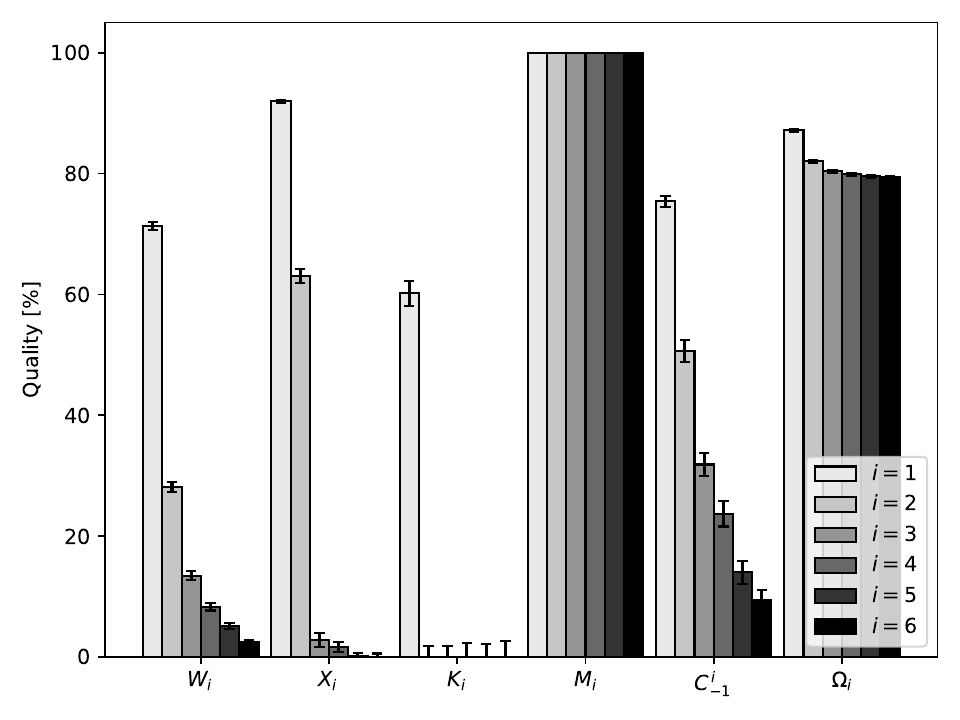}
		\caption{Low resolution.}
	\end{subfigure}%
	\begin{subfigure}{.5\textwidth}
		\centering
		\includegraphics[width=1\linewidth]{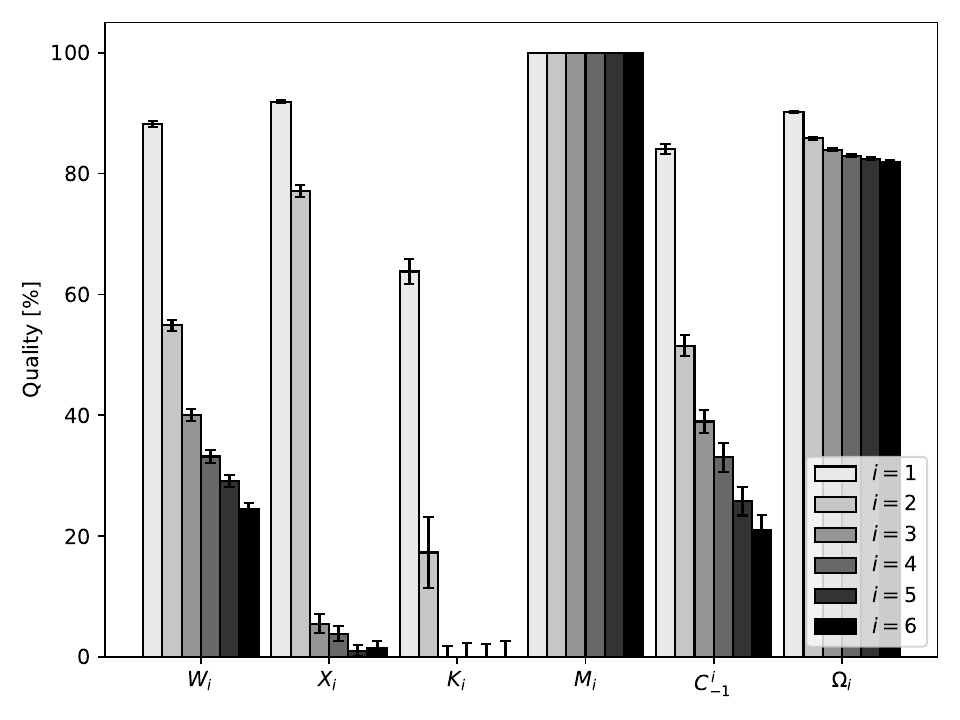}
		\caption{High resolution.}
	\end{subfigure}
	\caption{Quality criterion in percent of several sets and revolution periods. Note that $\Omega_i$ is used to refer to the quality criteria of $\Omega$ after $i$ periods.}
	\label{fig:quality}
\end{figure}

\Fig{fig:quality} represents the quality criterion for both low- and high-resolution mappings on the left and right, respectively. Each group of bars corresponds to a different collection of sets. For instance, \Xset{i} for $ i \in \llbracket1,6\rrbracket $. Differently, each color represents a different number of periods, from 1 (light) to 6 (dark). A bar's height represents the quality of that set in percent, while whiskers inform on the confidence interval at 95\% on such value. These estimations of the quality criterion are computed against the sample of size $ 10^{5}$ obtained with \gls*{GRATIS}.

The quality of acrobatic sets \Dset{i} is not shown since they are not mapped by the \gls*{DA} algorithm. In fact, their qualities are automatically set to 0\%. In addition, the quality of moon-crash sets \Mset{i} is always 100\%. Indeed, among the $ 10^{5} $ propagated trajectories, none of them crashed on Mars' moons. The quality of these empty sets is then $100\%$ by convention.

The quality criterion is always higher for the high-resolution mapping than for the low-resolution one. Especially for stable sets \Wset{i} and capture sets \Cset{i} since a larger part of these sets becomes inconsistent during the low-resolution propagation, as shown by \Fig{fig:classification_6_revolutions}. Furthermore, capture sets \Cset{i} are mapped in a better or equal way compared to stable sets \Wset{i}. This is due to the positioning of \Cset{i} farther from Mars \citep{merisio2021characterization}, than most of \Wset{i}. Thus, as highlighted in Figs.~\ref{fig:classification_2_revolutions} and~\ref{fig:classification_6_revolutions}, the consistency of this region is higher, which enables better quality. However, designing a specific mission with ballistic capture requires much more data than these two criteria and is beyond the scope of this work. Moreover, this work does not provide any additional property on the identified solutions, except that they belong to capture sets.

Except for those two groups of sets, only slight improvements are remarkable when switching from low resolution to high resolution. Indeed, apart from \Xset{1}, escape sets \Xset{i} and crash sets \Kset{i} fail to be accurately mapped, even though these regions of the search space are consistent. Indeed, the proposed classification algorithm does not track revolutions geometrically but uses period approximations. Consequently, the classification implies major approximations. Nevertheless, the poor mapping of these sets hardly impacts the quality of sets \Wset{i} due to their small sizes. Therefore, the global quality never drops below 80\% in low resolution. The overall quality can be improved at a high computational cost, compared to the point-wise mapping, thanks to high-resolution mapping. Such a mapping improves the global classification results by at least 4. These poor performances are mainly due to the inability of the \gls*{DA} mapping engine to track revolutions geometrically instead of temporally. Thus, \gls*{DA} mapping as presented in this work cannot outperform point-wise mapping in terms of classification precision, but it exceeds point-wise mapping concerning computational cost.

\section{Conclusion} \label{sec:conclusion}

This article presented a methodology to build cartography of \acrlong*{BC} sets using \acrlong*{DA} mapping. It proposes an alternative classification algorithm to sort the sub-domains produced by the \acrlong*{ADS} algorithm in newly defined capture sets. Moreover, instead of tracking revolutions with geometrical methods, as in \citet{luo2014constructing}, they are counted temporally, using regression on point-wise data. This work establishes a bridge between the point-wise mapping and \gls*{DA} mapping, allowing us to assess the performances of this new method. Hence, the introduction of two criteria. The first one is the consistency criterion, which represents the proportion of the search space where the mapping accuracy is guaranteed by the \gls*{ADS} algorithm. The second one is the quality criterion, which represents the success rate of the \gls*{DA} classification algorithm compared to the point-wise one.

Results show that \gls*{DA} mapping of \gls*{BC} sets could be performed on large search spaces, as \gls*{ADS} captured the dynamical variations on the whole domain. Furthermore, the consistency criterion shows that more than 87\% of the search space is guaranteed by the \gls*{ADS} algorithm as accurate, even for low-resolution mappings. Moreover, the quality criterion demonstrates that the global error rate of \gls*{DA} mapping is below 20\%. However, some small sets are either not or poorly mapped using revolution period regression, mainly due to the incapacity to count revolutions via geometrical information. This phenomenon does not disappear when the resolution of the mapping is increased. Nonetheless, \gls*{DA} mapping of \gls*{BC} sets can prove useful in various situations. For instance, it can be used for fast mapping of a plane, to then target a restricted domain of the search space. Either with a dense point-wise mapping or with a high-resolution \gls*{DA} mapping of the restricted area of interest.

\begin{dataavailability}

The data sets supporting the findings of this study are available on Zenodo\footnote{\url{https://zenodo.org/} [last accessed \lastdate]} with the identifiers \url{https://doi.org/10.5281/zenodo.6103720} and \url{https://doi.org/10.5281/zenodo.5940101}.

\end{dataavailability}

\begin{acknowledgements}

T.C. would like to thank the \gls*{DAER} of Politecnico di Milano for the host and the warm welcome, as well as the \gls*{TSAE} for the funding of his stay at Politecnico di Milano. G.M. and F.T. would like to acknowledge the \gls*{ERC} since part of this work has received funding from the \gls*{ERC} under the European Union’s Horizon 2020 research and innovation programme (grant agreement No.\,864697).

\end{acknowledgements}

\bibliography{references}    %

\end{document}